\documentclass[11pt]{article}
\usepackage{enumerate}
\usepackage{amssymb,a4wide,latexsym,makeidx,epsfig,fleqn}
\usepackage{amsthm}
\usepackage{amsmath}
\usepackage{enumerate}
\usepackage{graphicx}
\usepackage{float}
\usepackage{tikz}
\usepackage{hyperref}
\usetikzlibrary{positioning}
\usetikzlibrary{decorations.pathreplacing}
\usepackage{caption}
\allowdisplaybreaks[4]
\newtheorem{theorem}{Theorem}[section]

\newtheorem{lemma}[theorem]{Lemma}

\newtheorem{problem}[theorem]{Problem}

\def\qed{\hfill$\Box$\vspace{12pt}}

\begin{document}
\textwidth 150mm \textheight 225mm
\title{Tricyclic graphs with the second largest distance eigenvalue less than \( -\frac{1}{2}\)\thanks{Supported by National Natural Science Foundation of China (No. 12271439).} }
\author{{Kexin Yang$^{a,b,c}$, Ligong Wang$^{a,b,c}$\thanks{Corresponding author.
E-mail address: lgwangmath@163.com}}\\
{\small $^a$School of Mathematics and Statistics, Northwestern
Polytechnical University,}\\ {\small  Xi'an, Shaanxi 710129,
P.R. China}\\
{\small $^b$ Xi'an-Budapest Joint Research Center for Combinatorics, Northwestern
Polytechnical University,}\\
{\small Xi'an, Shaanxi 710129,
P.R. China}\\
{\small $^c$ MOE Key Laboratory for Complexity Science in Aerospace, Northwestern Polytechnical University,}\\
{\small Xi'an, Shaanxi 710129, P.R. China. }\\
{\small E-mail: yangkexi@mail.nwpu.edu.cn, lgwangmath@163.com} }
\date{}
\maketitle
\begin{center}
\begin{minipage}{120mm}
\vskip 0.3cm
\begin{center}
{\small {\bf Abstract}}
\end{center}
{\small 
	  Let $G$ be a simple connected graph with vertex set $V(G)=\{v_{1}, v_{2}, \ldots, v_{n}\}$. The distance $d_G(v_i,v_j)$ between two vertices $v_i$ and $v_j$ of $G$ is the length of a shortest path between $v_i$ and $v_j$. The distance matrix of $G$ is  defined as $D(G)=(d_G(v_i,v_j))_{n\times n}$. The second largest distance eigenvalue of \( G \) is the second largest eigenvalues of $D(G)$. Guo and Zhou [Discrete Math. 347(2024), 114082] proved that any connected graph with the second largest distance eigenvalue less than $-\frac{1}{2}$ is chordal, and characterize all bicyclic graphs and split graphs with the second largest distance eigenvalue less than $-\frac{1}{2}$. Based on this, we characterize all tricyclic graphs with the second largest distance eigenvalue less than $-\frac{1}{2}$.

}
\end{minipage}
\end{center}

\section{Introduction}
\label{introduction}

Throughout this paper, we consider simple connected graphs. Let $G$ be a graph with vertex set $V(G)=\{v_{1}, v_{2}, \ldots, v_{n}\}$ and edge set $E(G)=\{e_{1},e_{2}, \ldots, e_{m}\}$. We use $n=|V(G)|$ and $m=|E(G)|$ to denote the order and the size of $G$, respectively. If \( m = n + c - 1 \), then \( G \) is called a \( c \)-cyclic graph. If $c=0,1,2$ and 3, then $G$ is a tree, a unicyclic graph, a bicyclic graph and a tricyclic graph, respectively.
Let $P_n$, $C_n$, $K_n$ and $K_{t, n-t}$ be the path, the cycle, the complete graph and the complete bipartite graph of order $n$ respectively. 
 For a graph \( G \) with \( v \in V(G) \), \( N_G(v) \) is denoted the neighborhood of \( v \) in \( G \), and the degree of \( v \) in \( G \) is denoted by $ d_G(v) = |N_G(v)|$. Let $G[S]$ be the induced subgraph of $G$ whose vertex set is $S$ and edge set is consists of the edges which have both endpoints in $S$. A clique of a graph is a set of mutually adjacent vertices. An independent set of a graph is the set of vertices with no two of which are adjacent. A graph is a split graph when its vertex set can be partitioned into a clique and an independent set. A graph is chordal if every cycle of length at least four has a chord, where a chord is an edge joining two non-adjacent vertices of the cycle.
 
 The distance matrix $D(G)$ of a graph $G$ with vertex set $V(G)=\{v_{1}, v_{2}, \ldots, v_{n}\}$ is an $n\times n$ matrix whose $(i,j)$-entry is $d_G(v_i,v_j)$, where $d_G(v_i,v_j)$ is the distance between the vertices $v_i$ and $v_j$ in the graph $G$.
The eigenvalues of $D(G)$ are called the distance eigenvalues of $G$. In 1971, Graham and Pollak {\normalfont \cite{RF5}} firstly studied the distance eigenvalues of graphs. Their research revealed that the addressing problem in data communication systems is related to the number of negative distance eigenvalues of trees. In 1990, Merris {\normalfont \cite{RF11}} revealed a relationship between the distance eigenvalues and Laplacian eigenvalues of trees. These outstanding results made the spectral properties of distance matrix a research subject of great interest. Furthermore, growing attention has been given to the eigenvalues of distance matrix recently, but the emphasis primarily was on the first largest distance eigenvalue. Extensive findings regarding the distance spectra of graphs can see the excellent surveys {\normalfont \cite{RF2,RF10}}.

Since the distance matrix $D(G)$ of a graph $G$ with order $n$ is symmetric, the distance eigenvalues of  \( G \) are all real. We can order the eigenvalues of \( D(G) \) as \( \lambda_1(G) \ge\lambda_2(G) \ge \cdots \ge \lambda_n(G) \). Then \(\lambda_2(G) \) is called the second largest distance eigenvalue of graphs. Next, we present some results about the the second largest distance eigenvalue. In 1998, Fajtlowicz {\normalfont \cite{RF4}} conjectured that for any connected graph with independence number $\alpha _2\le 2$, $\lambda _2\left( G \right)$ is less than its number of triangles in $G$. Lin {\normalfont \cite{RF8}} comfirmed the conjecture and the equality holds if and only if $G\cong K_{2, 2}$. Xing and Zhou {\normalfont \cite{RF13}} provided rigorous lower bounds for the first and second largest distance eigenvalues of the $k$-th power of a connected graph. They determined all trees and unicyclic graphs $G$ such that $\lambda _2\left( G^2 \right) \le \frac{\sqrt{5}-3}{2}$. Besides, they also determined the unique $n$-vertex trees, the squares of which obtain minimum, second-minimum and third-minimum second largest distance eigenvalues in {\normalfont \cite{RF13}}. Liu et al. {\normalfont \cite{RF9}} proved that the graphs with $\lambda _2\left( G \right) \le \frac{17-\sqrt{329}}{2}$ can be determined by their distance spectra. Based on certain graph parameters, Alhevaz et al. {\normalfont \cite{RF1}} obtained some upper and lower bounds for the second largest eigenvalue of the generalized distance matrix of a graph. They also characterized the extremal graphs which attain these bounds. Xing and Zhou {\normalfont \cite{RF12}} gave a full characterization for connected graphs with $\lambda _2(G) \le -2+\sqrt{2}$. With the progress of the research on the second largest distance eigenvalue of a graph, a related problem comes into consideration.

\noindent\begin{problem}{\normalfont \cite{RF15}}
	Characterizing all connected graphs with the second largest distance eigenvalue at most
	\( -\frac{1}{2} \). 
\end{problem}
For the above problem, Xing and Zhou {\normalfont \cite{RF12}} determined  all  trees which satisfy $\lambda _2(G) \le -\frac{1}{2}$. Furthermore, they also characterized all unicyclic graphs with $\lambda _2(G) \le -\frac{1}{2}$, apart from a few exceptions. Xue et al. {\normalfont \cite{RF14}} determined all block graphs which satisfy $\lambda _2(G) \le -\frac{1}{2}$. Guo and Zhou \cite{RF6} characterized all bicyclic  graphs and split graphs with $\lambda _2(G)\le -\frac{1}{2}$, and proved that any connected  graph $G$ which satisfy $\lambda _2(G) \le  -\frac{1}{2}$ must be choral. 

 The rest of the paper is organized as follows. In Section 2, we give some important lemmas that will be used later. In Section 3, we characterize all tricyclic graphs whose second largest distance eigenvalue is less than \( -\frac{1}{2} \).

\section{Preliminaries}
\label{sec:Preliminaries}
In this section, we will present some important results that will be used in our subsequent arguments.

\begin{figure}[htbp]
	\hspace{1.2cm}
		\begin{minipage}[c]{0.3\textwidth}
			\begin{tikzpicture}[scale =1.5]
				\node[circle,fill=black,draw=black,inner sep=2.2pt] (v1) at (-0.4,0) {};
				\node[circle,fill=black,draw=black,inner sep=2.2pt] (v2) at (0.4,0) {};
				\node[circle,fill=black,draw=black,inner sep=2.2pt] (v5) at (1.4,0) {};
				\node[circle,fill=black,draw=black,inner sep=2.2pt] (v6) at (-1.4,0) {};
				\node[circle,fill=black,draw=black,inner sep=2.2pt] (v7) at (-1.6,0.8) {};
				\node[circle,fill=black,draw=black,inner sep=2.2pt] (v8) at (-2.2,0.4) {};
				\node[circle,fill=black,draw=black,inner sep=2.2pt] (v9) at (-2.1,-0.9) {};
				\node (v10) at (-0.9,-0.3) {};
				\node (v11) at (0.9,-0.3) {};
				\node (v12) at (0,-1.0) {};
				\node[below] at (v1) { };
				\node[left] at (v2) { };
				\node[left] at (v5) { };
				\node[below] at (v6) { };
				\node[right] at (v7) { };
				\node[right] at (v8) { };
				\node[right] at (v9) { };
				\node[below] at (v10) {\large$K_p$};	
				\node[below] at (v11) {\large$K_q$};	
				\node[below] at (v12) {\large$BG(p,q,3,2,2)$};	
				\draw [line width=1.2pt](v1) -- (v2);
				\draw[line width=1.2pt] (v6) -- (v7);
				\draw [line width=1.2pt](v6) -- (v8);
				\draw [line width=1.2pt](v6) -- (v9);
				\draw [line width=1.2pt](v7) -- (v8);
				\draw [line width=1.2pt](-0.9,0) ellipse (0.65cm and 0.17cm);
				\draw [line width=1.2pt](0.9,0) ellipse (0.65cm and 0.17cm);
				\fill (-0.7,0) circle (1pt); 
				\fill (-0.9,0) circle (1pt); 
				\fill (-1.1,0) circle (1pt); 
				\fill (0.7,0) circle (1pt); 
				\fill (0.9,0) circle (1pt); 
				\fill (1.1,0) circle (1pt); 
			\end{tikzpicture}
		\end{minipage}
		\hspace{0.2\textwidth}
		\begin{minipage}[c]{0.3\textwidth}
			\begin{tikzpicture}[scale =1.5]
				\node[circle,fill=black,draw=black,inner sep=2.2pt] (v1) at (0,0) {};
				\node[circle,fill=black,draw=black,inner sep=2.2pt] (v2) at (-0.9,0) {};
				\node[circle,fill=black,draw=black,inner sep=2.2pt] (v3) at (0.9,0) {};
				\node[circle,fill=black,draw=black,inner sep=2.2pt] (v4) at (1.4,-0.9) {};
				\node[circle,fill=black,draw=black,inner sep=2.2pt] (v5) at (-1.4,-0.9) {};
				\node[circle,fill=black,draw=black,inner sep=2.2pt] (v6) at (-1.1,0.8) {};
				\node[circle,fill=black,draw=black,inner sep=2.2pt] (v7) at (1.1,0.8) {};
				\node[circle,fill=black,draw=black,inner sep=2.2pt] (v8) at (1.6,0.3) {};
				\node[circle,fill=black,draw=black,inner sep=2.2pt] (v9) at (-1.6,0.3) {};
				\node (v10) at (0,-0.95) {};
				\node[below] at (v1) { };
				\node[left] at (v2) { };
				\node[left] at (v3) { };
				\node[left] at (v4) { };
				\node[left] at (v5) { };
				\node[below] at (v6) { };
				\node[right] at (v7) { };
				\node[right] at (v8) { };
				\node[right] at (v9) { };
				\node[below] at (v10) {\large$BGA$};	
				\draw [line width=1.2pt](v1) -- (v2);
				\draw[line width=1.2pt] (v1) -- (v3);
				\draw [line width=1.2pt](v3) -- (v7);
				\draw [line width=1.2pt](v3) -- (v8);
				\draw [line width=1.2pt](v3) -- (v4);
				\draw[line width=1.2pt] (v2) -- (v6);
				\draw [line width=1.2pt](v2) -- (v9);
				\draw [line width=1.2pt](v2) -- (v5) ;
				\draw [line width=1.2pt](v9) -- (v6);
				\draw [line width=1.2pt](v7) -- (v8);	
			\end{tikzpicture}
		\end{minipage}
	\caption{Graphs $BG(p,q,3,2,2)$ and $BGA$ } 
	\label{figure:bga} 
\end{figure}
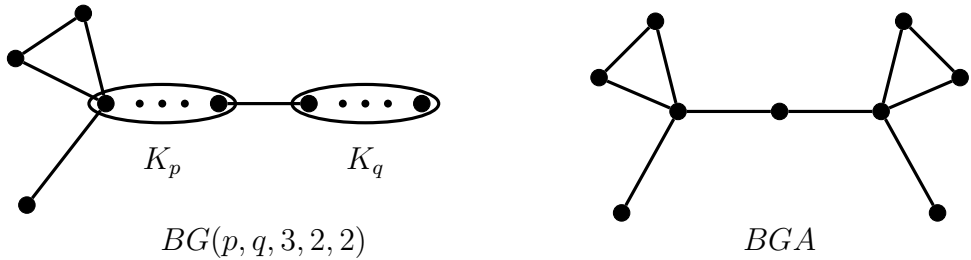

A block of a graph is a maximal connected subgraph without cut vertex. If all blocks of a connected graph are clique, then we call \( G \) a block graph (or clique tree). A block star is a block graph whose all blocks share a common vertex. A block graph \( G \) is loose if for every vertex $v\in V(G)$, the vertex $v$ is contained in at most two blocks of $G$. Let $BG(p,q,3,2,2)$, $p,q\ge2$, and \( BGA \) be two block graphs as depicted in Figure. \ref{figure:bga}.

\noindent\begin{lemma}\label{th:ch-1}{\normalfont(\cite{RF12})}
	Let $G$ be a block graph with the second largest distance eigenvalue $\lambda_2(G)$. Then $\lambda_2(G) < -\frac{1}{2}$ if and only if
	\begin{itemize}
		\item $G$ is a block star, or
		\item $G$ is a loose block graph, or
		\item $G$ is a nontrivial connected induced subgraph of $BG(p,q,3,2,2)$, where $p,q\ge2$, or
		\item $G$ is a nontrivial connected induced subgraph of $BGA$.
	\end{itemize}
\end{lemma}

For a real symmetric matrix $M$ of order $n$, then the eigenvalues of $M$ can be ordered as $\rho _1(M)\ge \rho _2(M)\ge \cdots \ge \rho _n(M)$.
\noindent\begin{lemma}\label{th:ch-2}{\normalfont(\cite{RF7})}
	(Cauchy's interlacing theorem). Let $A$ be a symmetric matrix of order $n$. If $B$ is a principal submatrix of $A$ of order $m$, then $\rho_i(A) \geq \rho_i(B) \geq \rho_{n-m+i}(A)$ for $1 \leq i \leq m$.
\end{lemma}

Let $H$ be a distance-preserving subgraph of a graph $G$ when $H$ is a connected induced subgraph of $G$ such that $d_H(u, v) = d_G(u, v)$ for  all $\{u, v\} \subseteq V(H)$. In this situation, $D(H)$ is a principal submatrix of $D(G)$, by Lemma \ref{th:ch-2}, we have $\lambda_2(G) \geq \lambda_2(H)$. To be specific, if $G$ is a connected graph which satisfies $\lambda_2(G) < -\frac{1}{2}$, then for all distance-preserving subgraphs $H$ of $G$ which satisfy $\lambda_2(H) < -\frac{1}{2}$ . 

In this paper, we will characterize tricyclic graphs with $\lambda_2(G) < -\frac{1}{2}$. then we name a graph $H$ a forbidden subgraph of $G$ when $H$ is a distance-preserving  subgraph of $G$ but $\lambda_2(H) \geq -\frac{1}{2}$. We show in Figure. \ref{fig:forbidden} the set of forbidden subgraphs $F_1$ to $ F_{13}$, which are required in the proof, and the graphs are accompanied by their second largest distance eigenvalue information directly beneath them.

Let $\mathcal{\varPi} = \{V_1,V_2,\ldots,V_m\}$ be a partition of $\{1,2,\ldots, n\}$.
If the element $b_{ij}$ of the matrix $B_\varPi$ is the average row sum of the block $A_{i,j}$ of a real symmetric matrix $A$ with the rows and columns indexed by $V_i$ and $V_j$, respectively, then $B_\varPi$ is a quotient matrix of $A$ with respect to $\mathcal{\varPi}$. If each block of the partition has constant row sums, then $\mathcal{\varPi}$ is called an equitable partition of $A$.

\noindent\begin{lemma}\label{th:ch-3}{\normalfont(\cite{RF3})} Let $A$ be an $n\times n$  real symmetric matrix and let $\varPi$ be a partition of $\{1,2,\ldots, n\}$. If the corresponding partition of $A$ is equitable and its quotient matrix is $B_{\varPi}$, then the spectrum of $B_{\varPi}$ is a sub(multi) set of the spectrum of $A$.
\end{lemma}

\noindent\begin{lemma}\label{th:ch-4}{\normalfont(\cite{RF7})}
	If $G$ is a connected graph with \( \lambda_2(G) \) $<$ \(-\frac{1}{2}\), then $G$ is chordal.
\end{lemma}
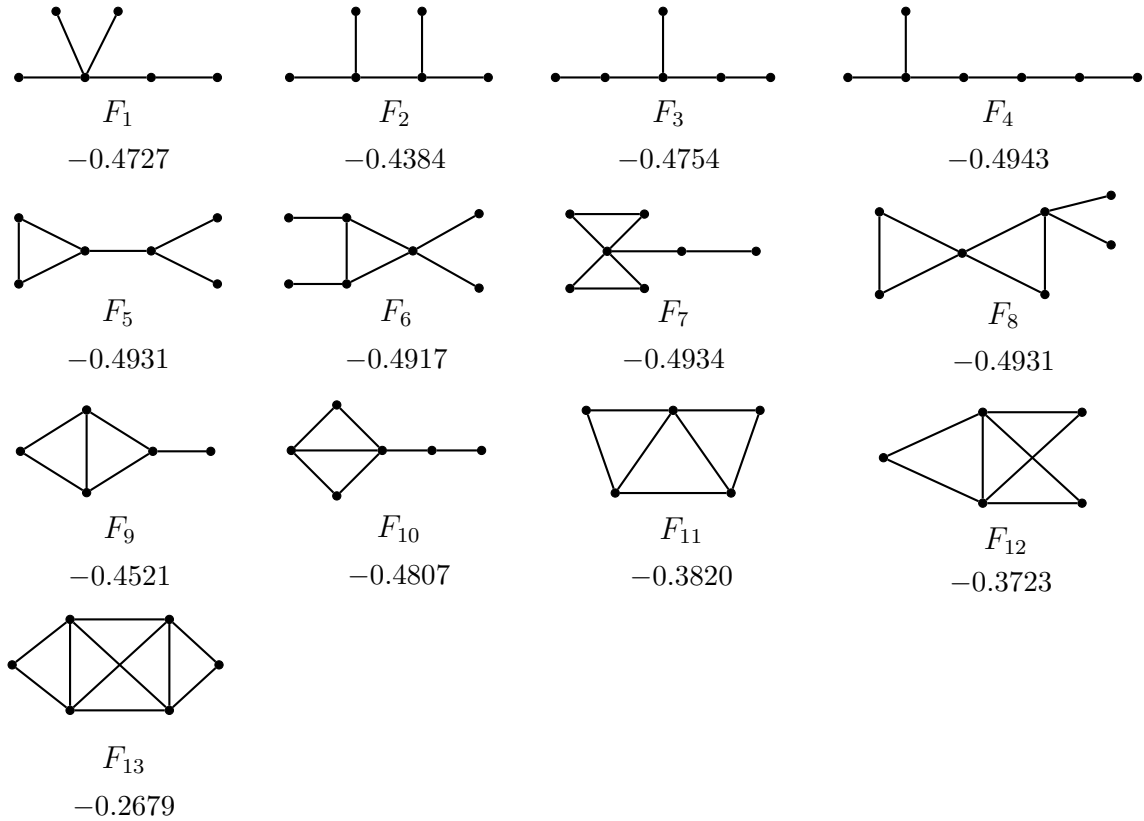
\begin{figure}[h]
	\begin{tikzpicture}
		\node (fig1) {	  	
			\begin{tikzpicture}[scale=1.1]
				\node[circle,fill=black,draw=black,inner sep=1.1pt] (v1) at (-0.4,0) {};
				\node[circle,fill=black,draw=black,inner sep=1.1pt] (v2) at (0.4,0) {};
				\node[circle,fill=black,draw=black,inner sep=1.1pt] (v3) at (1.2,0) {};
				\node[circle,fill=black,draw=black,inner sep=1.1pt] (v4) at (-1.2,0) {};
				\node[circle,fill=black,draw=black,inner sep=1.1pt] (v5) at (0,0.8) {};
				\node[circle,fill=black,draw=black,inner sep=1.1pt] (v6) at (-0.75,0.8) {};
				\node (v7) at (0,-0.15) {};
				\node (v8) at (0,-0.75) {};
				\node[below] at (v1) { };
				\node[left] at (v2) { };
				\node[left] at (v3) { };
				\node[below] at (v4) { };
				\node[right] at (v5) { };
				\node[right] at (v6) { };
				\node[below] at (v7) {\large$F_1$};	
				\node[below] at (v8) {$-0.4727$};	
				
				\draw [line width=0.8pt] (v1) -- (v2);
				\draw [line width=0.8pt] (v2) -- (v3);
				\draw [line width=0.8pt] (v1) -- (v4);
				\draw [line width=0.8pt] (v1) -- (v6);
				\draw [line width=0.8pt] (v1) -- (v5);
				
			\end{tikzpicture}
		};
		\node[right=0.5cm of fig1] (fig2) {
			\begin{tikzpicture}[scale =1.1]
				\node[circle,fill=black,draw=black,inner sep=1.1pt] (v1) at (-0.4,0) {};
				\node[circle,fill=black,draw=black,inner sep=1.1pt] (v2) at (0.4,0) {};
				\node[circle,fill=black,draw=black,inner sep=1.1pt] (v3) at (1.2,0) {};
				\node[circle,fill=black,draw=black,inner sep=1.1pt] (v4) at (-1.2,0) {};
				\node[circle,fill=black,draw=black,inner sep=1.1pt] (v5) at (-0.4,0.8) {};
				\node[circle,fill=black,draw=black,inner sep=1.1pt] (v6) at (0.4,0.8) {};
				\node (v7) at (0,-0.15) {};
				\node (v8) at (0,-0.75) {};
				\node[below] at (v1) { };
				\node[left] at (v2) { };
				\node[left] at (v3) { };
				\node[below] at (v4) { };
				\node[right] at (v5) { };
				\node[right] at (v6) { };
				\node[below] at (v7) {\large$F_2$};	
				\node[below] at (v8) {$-0.4384$};	
				
				\draw [line width=0.8pt](v1) -- (v2);
				\draw[line width=0.8pt] (v2) -- (v3);
				\draw [line width=0.8pt](v1) -- (v4);
				\draw [line width=0.8pt](v1) -- (v5);
				\draw [line width=0.8pt](v2) -- (v6);
				
			\end{tikzpicture}
		};
		\node[right=0.5cm of fig2] (fig3) {
			\begin{tikzpicture}[scale =1.1]
				\node[circle,fill=black,draw=black,inner sep=1.1pt] (v1) at (0,0) {};
				\node[circle,fill=black,draw=black,inner sep=1.1pt] (v2) at (0.7,0) {};
				\node[circle,fill=black,draw=black,inner sep=1.1pt] (v3) at (-0.7,0) {};
				\node[circle,fill=black,draw=black,inner sep=1.1pt] (v4) at (1.3,0) {};
				\node[circle,fill=black,draw=black,inner sep=1.1pt] (v5) at (-1.3,0) {};
				\node[circle,fill=black,draw=black,inner sep=1.1pt] (v6) at (0,0.8) {};
				\node (v7) at (0,-0.15) {};
				\node (v8) at (0,-0.75) {};
				\node[below] at (v1) { };
				\node[left] at (v2) { };
				\node[left] at (v3) { };
				\node[below] at (v4) { };
				\node[right] at (v5) { };
				\node[right] at (v6) { };
				\node[below] at (v7) {\large$F_3$};	
				\node[below] at (v8) {$-0.4754$};	
				
				\draw [line width=0.8pt](v1) -- (v6);
				\draw[line width=0.8pt] (v1) -- (v2);
				\draw [line width=0.8pt](v1) -- (v3);
				\draw [line width=0.8pt](v2) -- (v4);
				\draw [line width=0.8pt](v3) -- (v5);
				
			\end{tikzpicture}
		};
		\node[right=0.5cm of fig3] (fig4) {
			\begin{tikzpicture}[scale =1.1]
				\node[circle,fill=black,draw=black,inner sep=1.1pt] (v0) at (-0.35,0) {};
				\node[circle,fill=black,draw=black,inner sep=1.1pt] (v1) at (0.35,0) {};
				\node[circle,fill=black,draw=black,inner sep=1.1pt] (v2) at (-1.05,0) {};
				\node[circle,fill=black,draw=black,inner sep=1.1pt] (v3) at (1.05,0) {};
				\node[circle,fill=black,draw=black,inner sep=1.1pt] (v4) at (-1.75,0) {};
				\node[circle,fill=black,draw=black,inner sep=1.1pt] (v5) at (1.75,0) {};
				\node[circle,fill=black,draw=black,inner sep=1.1pt] (v6) at (-1.05,0.8) {};
				\node (v7) at (0,-0.15) {};
				\node (v8) at (0,-0.75) {};
				\node[below] at (v0) { };
				\node[below] at (v1) { };
				\node[left] at (v2) { };
				\node[left] at (v3) { };
				\node[below] at (v4) { };
				\node[right] at (v5) { };
				\node[right] at (v6) { };
				\node[below] at (v7) {\large$F_4$};	
				\node[below] at (v8) {$-0.4943$};	
				
				\draw [line width=0.8pt](v0) -- (v1);
				\draw[line width=0.8pt] (v1) -- (v3);
				\draw [line width=0.8pt](v3) -- (v5);
				\draw [line width=0.8pt](v0) -- (v2);
				\draw [line width=0.8pt](v2) -- (v6);
				\draw [line width=0.8pt](v2) -- (v4);
			\end{tikzpicture}
		};
		\node[below=0.1cm of fig1] (fig5) {
			\begin{tikzpicture}[scale =1.1]
				\node[circle,fill=black,draw=black,inner sep=1.1pt] (v0) at (-0.4,0) {};
				\node[circle,fill=black,draw=black,inner sep=1.1pt] (v1) at (0.4,0) {};
				\node[circle,fill=black,draw=black,inner sep=1.1pt] (v2) at (1.2,0.4) {};
				\node[circle,fill=black,draw=black,inner sep=1.1pt] (v3) at (1.2,-0.4) {};
				\node[circle,fill=black,draw=black,inner sep=1.1pt] (v4) at (-1.2,0.4) {};
				\node[circle,fill=black,draw=black,inner sep=1.1pt] (v5) at (-1.2,-0.4) {};
				\node (v7) at (0,-0.4) {};
				\node (v8) at (0,-1) {};
				\node[below] at (v0) { };
				\node[below] at (v1) { };
				\node[left] at (v2) { };
				\node[left] at (v3) { };
				\node[below] at (v4) { };
				\node[right] at (v5) { };
				\node[below] at (v7) {\large$F_5$};	
				\node[below] at (v8) {$-0.4931$};	
				
				\draw [line width=0.8pt](v0) -- (v1);
				\draw[line width=0.8pt] (v1) -- (v3);
				\draw [line width=0.8pt](v1) -- (v2);
				\draw [line width=0.8pt](v0) -- (v4);
				\draw [line width=0.8pt](v0) -- (v5);
				\draw [line width=0.8pt](v4) -- (v5);
			\end{tikzpicture}
		};
		\node[below=0.1cm of fig2] (fig6) {
			\begin{tikzpicture}[scale =1.1]
				\node[circle,fill=black,draw=black,inner sep=1.1pt] (v0) at (-0.1,0) {};
				\node[circle,fill=black,draw=black,inner sep=1.1pt] (v1) at (0.7,0.45) {};
				\node[circle,fill=black,draw=black,inner sep=1.1pt] (v2) at (0.7,-0.45) {};
				\node[circle,fill=black,draw=black,inner sep=1.1pt] (v3) at (-0.9,0.4) {};
				\node[circle,fill=black,draw=black,inner sep=1.1pt] (v4) at (-0.9,-0.4) {};
				\node[circle,fill=black,draw=black,inner sep=1.1pt] (v5) at (-1.6,0.4) {};
				\node[circle,fill=black,draw=black,inner sep=1.1pt] (v6) at (-1.6,-0.4) {};
				\node (v7) at (-0.3,-0.4) {};
				\node (v8) at (-0.3,-1) {};
				\node[below] at (v0) { };
				\node[below] at (v1) { };
				\node[left] at (v2) { };
				\node[left] at (v3) { };
				\node[below] at (v4) { };
				\node[right] at (v5) { };
				\node[right] at (v6) { };
				\node[below] at (v7) {\large$F_6$};	
				\node[below] at (v8) {$-0.4917$};	
				
				\draw [line width=0.8pt](v0) -- (v1);
				\draw[line  width=0.8pt] (v0) -- (v2);
				\draw [line width=0.8pt](v0) -- (v3);
				\draw [line width=0.8pt](v0) -- (v4);
				\draw [line width=0.8pt](v3) -- (v5);
				\draw [line width=0.8pt](v4) -- (v6);
				\draw [line width=0.8pt](v3) -- (v4);
			\end{tikzpicture}
		};
		\node[below=0.05cm of fig3] (fig7) {
			\begin{tikzpicture}[scale =1.1]
				\node[circle,fill=black,draw=black,inner sep=1.1pt] (v0) at (0.4,0) {};
				\node[circle,fill=black,draw=black,inner sep=1.1pt] (v1) at (1.3,0) {};
				\node[circle,fill=black,draw=black,inner sep=1.1pt] (v2) at (-0.5,0) {};
				\node[circle,fill=black,draw=black,inner sep=1.1pt] (v3) at (-0.05,0.45) {};
				\node[circle,fill=black,draw=black,inner sep=1.1pt] (v4) at (-0.05,-0.45) {};
				\node[circle,fill=black,draw=black,inner sep=1.1pt] (v5) at (-0.95,0.45) {};
				\node[circle,fill=black,draw=black,inner sep=1.1pt] (v6) at (-0.95,-0.45) {};
				\node (v7) at (0.3,-0.4) {};
				\node (v8) at (0.3,-1) {};
				\node[below] at (v0) { };
				\node[below] at (v1) { };
				\node[left] at (v2) { };
				\node[left] at (v3) { };
				\node[below] at (v4) { };
				\node[right] at (v5) { };
				\node[right] at (v6) { };
				\node[below] at (v7) {\large$F_7$};	
				\node[below] at (v8) {$-0.4934$};	
				
				\draw [line width=0.8pt](v0) -- (v1);
				\draw[line  width=0.8pt] (v0) -- (v2);
				\draw [line width=0.8pt](v2) -- (v3);
				\draw [line width=0.8pt](v2) -- (v4);
				\draw [line width=0.8pt](v2) -- (v5);
				\draw [line width=0.8pt](v2) -- (v6);
				\draw [line width=0.8pt](v3) -- (v5);
				\draw [line width=0.8pt](v4) -- (v6);
			\end{tikzpicture}
		};
		\node[below=-0.2cm of fig4] (fig8) {
			\begin{tikzpicture}[scale =1.1]
				\node[circle,fill=black,draw=black,inner sep=1.1pt] (v0) at (-0.2,0) {};
				\node[circle,fill=black,draw=black,inner sep=1.1pt] (v1) at (0.8,0.5) {};
				\node[circle,fill=black,draw=black,inner sep=1.1pt] (v2) at (0.8,-0.5) {};
				\node[circle,fill=black,draw=black,inner sep=1.1pt] (v3) at (1.6,0.7) {};
				\node[circle,fill=black,draw=black,inner sep=1.1pt] (v4) at (1.6,0.1) {};
				\node[circle,fill=black,draw=black,inner sep=1.1pt] (v5) at (-1.2,0.5) {};
				\node[circle,fill=black,draw=black,inner sep=1.1pt] (v6) at (-1.2,-0.5) {};
				\node (v7) at (0.3,-0.4) {};
				\node (v8) at (0.3,-1) {};
				\node[below] at (v0) { };
				\node[below] at (v1) { };
				\node[left] at (v2) { };
				\node[left] at (v3) { };
				\node[below] at (v4) { };
				\node[right] at (v5) { };
				\node[right] at (v6) { };
				\node[below] at (v7) {\large$F_8$};	
				\node[below] at (v8) {$-0.4931$};	
				
				\draw [line width=0.8pt](v0) -- (v1);
				\draw[line  width=0.8pt] (v0) -- (v2);
				\draw [line width=0.8pt](v0) -- (v5);
				\draw [line width=0.8pt](v0) -- (v6);
				\draw [line width=0.8pt](v1) -- (v3);
				\draw [line width=0.8pt](v1) -- (v4);
				\draw [line width=0.8pt](v1) -- (v2);
				\draw [line width=0.8pt](v5) -- (v6);
			\end{tikzpicture}
		};
		\node[below=0cm of fig5] (fig9) {
			\begin{tikzpicture}[scale =1.1]
				\node[circle,fill=black,draw=black,inner sep=1.1pt] (v1) at (0.4,0) {};
				\node[circle,fill=black,draw=black,inner sep=1.1pt] (v2) at (-1.2,0) {};
				\node[circle,fill=black,draw=black,inner sep=1.1pt] (v3) at (-0.4,0.5) {};
				\node[circle,fill=black,draw=black,inner sep=1.1pt] (v4) at (-0.4,-0.5) {};
				\node[circle,fill=black,draw=black,inner sep=1.1pt] (v5) at (1.1,0) {};
				\node (v7) at (0,-0.6) {};
				\node (v8) at (0,-1.2) {};
				\node[below] at (v1) { };
				\node[left] at (v2) { };
				\node[left] at (v3) { };
				\node[below] at (v4) { };
				\node[right] at (v5) { };
				\node[below] at (v7) {\large$F_9$};	
				\node[below] at (v8) {$-0.4521$};	
				
				\draw [line width=0.8pt](v1) -- (v5);
				\draw[line  width=0.8pt] (v1) -- (v3);
				\draw [line width=0.8pt](v1) -- (v4);
				\draw [line width=0.8pt](v2) -- (v3);
				\draw [line width=0.8pt](v2) -- (v4);
				\draw [line width=0.8pt](v3) -- (v4);
			\end{tikzpicture}
		};
		\node[below=0cm of fig6] (fig10) {
			\begin{tikzpicture}[scale =1.1]
				\node[circle,fill=black,draw=black,inner sep=1.1pt] (v1) at (0.2,0) {};
				\node[circle,fill=black,draw=black,inner sep=1.1pt] (v2) at (0.8,0) {};
				\node[circle,fill=black,draw=black,inner sep=1.1pt] (v3) at (1.4,0) {};
				\node[circle,fill=black,draw=black,inner sep=1.1pt] (v4) at (-0.35,0.55) {};
				\node[circle,fill=black,draw=black,inner sep=1.1pt] (v5) at (-0.35,-0.55) {};
				\node[circle,fill=black,draw=black,inner sep=1.1pt] (v6) at (-0.9,0) {};
				\node (v7) at (0.4,-0.6) {};
				\node (v8) at (0.4,-1.2) {};
				\node[below] at (v1) { };
				\node[left] at (v2) { };
				\node[left] at (v3) { };
				\node[below] at (v4) { };
				\node[right] at (v5) { };
				\node[below] at (v7) {\large$F_{10}$};	
				\node[below] at (v8) {$-0.4807$};	
				
				\draw [line width=0.8pt](v1) -- (v2);
				\draw[line  width=0.8pt](v1) -- (v5);
				\draw [line width=0.8pt](v1) -- (v4);
				\draw [line width=0.8pt](v1) -- (v6);
				\draw [line width=0.8pt](v2) -- (v3);
				\draw [line width=0.8pt](v4) -- (v6);
				\draw [line width=0.8pt](v5) -- (v6);
			\end{tikzpicture}
		};
		\node[below=0cm of fig7] (fig11) {
			\begin{tikzpicture}[scale =1.1]
				\node[circle,fill=black,draw=black,inner sep=1.1pt] (v1) at (0,0) {};
				\node[circle,fill=black,draw=black,inner sep=1.1pt] (v2) at (1.05,0) {};
				\node[circle,fill=black,draw=black,inner sep=1.1pt] (v3) at (-1.05,0) {};
				\node[circle,fill=black,draw=black,inner sep=1.1pt] (v4) at (-0.7,-1) {};
				\node[circle,fill=black,draw=black,inner sep=1.1pt] (v5) at (0.7,-1) {};
				\node (v7) at (0.1,-1.1) {};
				\node (v8) at (0.1,-1.7) {};
				\node[below] at (v1) { };
				\node[left] at (v2) { };
				\node[left] at (v3) { };
				\node[below] at (v4) { };
				\node[right] at (v5) { };
				\node[below] at (v7) {\large$F_{11}$};	
				\node[below] at (v8) {$-0.3820$};	
				
				\draw [line width=0.8pt](v1) -- (v2);
				\draw[line  width=0.8pt](v1) -- (v3);
				\draw[line  width=0.8pt](v1) -- (v5);
				\draw [line width=0.8pt](v1) -- (v4);
				\draw [line width=0.8pt](v2) -- (v5);
				\draw [line width=0.8pt](v4) -- (v5);
				\draw [line width=0.8pt](v4) -- (v3);
			\end{tikzpicture}
		};
		\node[below=0cm of fig8] (fig12) {
			\begin{tikzpicture}[scale =1.1]
				\hspace{-0.2cm}
				\node[circle,fill=black,draw=black,inner sep=1.1pt] (v1) at (-0.8,0.65) {};
				\node[circle,fill=black,draw=black,inner sep=1.1pt] (v2) at (0.4,0.65) {};
				\node[circle,fill=black,draw=black,inner sep=1.1pt] (v3) at (0.4,-0.45) {};
				\node[circle,fill=black,draw=black,inner sep=1.1pt] (v4) at (-0.8,-0.45) {};
				\node[circle,fill=black,draw=black,inner sep=1.1pt] (v5) at (-2,0.1) {};
				\node (v7) at (-0.52,-0.55) {};
				\node (v8) at (-0.58,-1.1) {};
				\node[below] at (v1) { };
				\node[left] at (v2) { };
				\node[left] at (v3) { };
				\node[below] at (v4) { };
				\node[right] at (v5) { };
				\node[below] at (v7) {\large$F_{12}$};	
				\node[below] at (v8) {$-0.3723$};	
				
				\draw [line width=0.8pt](v1) -- (v2);
				\draw[line  width=0.8pt](v1) -- (v3);
				\draw[line  width=0.8pt](v1) -- (v5);
				\draw [line width=0.8pt](v1) -- (v4);
				\draw [line width=0.8pt](v2) -- (v4);
				\draw [line width=0.8pt](v4) -- (v5);
				\draw [line width=0.8pt](v4) -- (v3);
			\end{tikzpicture}
		};
		\node[below=-0.1cm of fig9] (fig13) {
			\begin{tikzpicture}[scale =1.1]
				\node[circle,fill=black,draw=black,inner sep=1.1pt] (v1) at (-0.6,0.55) {};
				\node[circle,fill=black,draw=black,inner sep=1.1pt] (v2) at (0.6,0.55) {};
				\node[circle,fill=black,draw=black,inner sep=1.1pt] (v3) at (0.6,-0.55) {};
				\node[circle,fill=black,draw=black,inner sep=1.1pt] (v4) at (-0.6,-0.55) {};
				\node[circle,fill=black,draw=black,inner sep=1.1pt] (v5) at (-1.3,0) {};
				\node[circle,fill=black,draw=black,inner sep=1.1pt] (v6) at (1.2,0) {};
				\node (v7) at (0.05,-0.8) {};
				\node (v8) at (0.05,-1.4) {};
				\node[below] at (v1) { };
				\node[left] at (v2) { };
				\node[left] at (v3) { };
				\node[below] at (v4) { };
				\node[right] at (v5) { };
				\node[below] at (v7) {\large$F_{13}$};	
				\node[below] at (v8) {$-0.2679$};	
				
				\draw [line width=0.8pt](v1) -- (v2);
				\draw[line  width=0.8pt](v1) -- (v3);
				\draw[line  width=0.8pt](v1) -- (v5);
				\draw [line width=0.8pt](v1) -- (v4);
				\draw [line width=0.8pt](v2) -- (v4);
				\draw [line width=0.8pt](v4) -- (v5);
				\draw [line width=0.8pt](v4) -- (v3);
				\draw [line width=0.8pt](v2) -- (v6);
				\draw [line width=0.8pt](v3) -- (v6);
				\draw [line width=0.8pt](v2) -- (v3);
			\end{tikzpicture}
		};
	\end{tikzpicture}
	\caption{Forbidden subgraphs $F_i$ with $\lambda_2 \geq -\frac{1}{2}$ ($i$=1,...,13) } 
	\label{fig:forbidden} 
\end{figure}
Let $G$ be a tricyclic graph, we use $B(G)$ to denote the base of $G$ which is the minimal tricyclic subgraph of $G$. Obviously, $B(G)$ is the unique tricyclic subgraph of $G$ that contains no pendant vertices, and $G$ can be obtained from $B(G)$ by attaching trees to some vertices of $B(G)$.
\noindent\begin{lemma}\label{th:ch-5}{\normalfont(\cite{RF15})}
	If $G$ is a tricyclic graph, then $G$ have the following four types of bases :$G_j^3$ ($j$=1,...,7), $G_j^4$ ($j$=1,...,4), $G_j^6$ ($j$=1,...,3), $G_1^7$, as shown in Figure.\ref{figure:G3}.
\end{lemma}
\begin{figure}[htbp]
	\begin{tikzpicture}
		\node (fig1) {	  	
			\begin{tikzpicture}[scale =1]
				\hspace{0.1cm}
				\draw[line width=1.1pt, thick, black] (0,0) circle (0.5);
				\draw[line width=1.1pt, thick, black] (-1,0) circle (0.5);
				\draw[line width=1.1pt, thick, black] (1,0) circle (0.5); 
				\node[circle,fill=black,draw=black,inner sep=1.5pt] (v1) at (0.5,0) {};
				\node[circle,fill=black,draw=black,inner sep=1.5pt] (v2) at (-0.5,0) {};
				\node at (0,-1.2) {\large$G_1^3$};
				
			\end{tikzpicture}
		};
		\node[right=0.6cm of fig1] (fig2) {
			\begin{tikzpicture}[scale =1]
				\hspace{0.6cm}
				\draw[line width=1.1pt, thick, black] (-0.5,0) circle (0.5);
				\draw[line width=1.1pt, thick, black] (-1.5,0) circle (0.5);
				\draw[line width=1.1pt, thick, black] (1.5,0) circle (0.5); 
				\node[circle,fill=black,draw=black,inner sep=1.5pt] (v1) at (-0.085,0) {};
				\node[circle,fill=black,draw=black,inner sep=1.5pt] (v2) at (-1.08,0) {};
				\node[circle,fill=black,draw=black,inner sep=1.5pt] (v3) at (0.92,0) {};
				\draw [dash pattern=on 3pt off 2pt, line width=1.1pt](v1) -- (v3);
				\node at (0.25,0.3) {$P$};
				\node at (0.25,-1.2) {\large$G_2^3$};

			\end{tikzpicture}
		};
		\node[right=0.35cm of fig2] (fig3) {
			\begin{tikzpicture}[scale =1]
				\hspace{1.3cm}
				\draw[line width=1.1pt, thick, black] (0,0) circle (0.5);
				\draw[line width=1.1pt, thick, black] (2,0) circle (0.5);
				\draw[line width=1.1pt, thick, black] (-2,0) circle (0.5); 
				\node[circle,fill=black,draw=black,inner sep=1.5pt] (v1) at (0.41,0) {};
				\node[circle,fill=black,draw=black,inner sep=1.5pt] (v2) at (-0.58,0) {};
				\node[circle,fill=black,draw=black,inner sep=1.5pt] (v3) at (1.42,0) {};
				\node[circle,fill=black,draw=black,inner sep=1.5pt] (v4) at (-1.58,0) {};
				\draw [dash pattern=on 3pt off 2pt, line width=1.1pt](v1) -- (v3);
				\draw [dash pattern=on 3pt off 2pt, line width=1.1pt](v2) -- (v4);
				\node at (0.7,0.3) {$P^2$};
				\node at (-1.3,0.3) {$P^1$};
				\node at (-0.2,-1.2) {\large$G_3^3$};
			\end{tikzpicture}
		};
		
		\node[below=0.4cm of fig1] (fig4) {
			\begin{tikzpicture}[scale =1]
				\def\R{1.2}                
				\def\k{3}                  
				\draw[line width=1.1pt,black,domain=0:360,samples=500,smooth] 
				plot ({\x}:{\R*cos(\k*\x)});
				\node at (0,-1.5) {\large$G_4^3$};
				\node[circle,fill=black,draw=black,inner sep=1.5pt]  at (0,0.1) {};
			\end{tikzpicture}
		};
		\node[right=0.5cm of fig4] (fig5) {
			\begin{tikzpicture}[scale =1]
				\draw[line width=1.1pt, thick, black] (0,-0.5) circle (0.5);
				\draw[line width=1.1pt, thick, black] (0,0.5) circle (0.5);
				\draw[line width=1.1pt, thick, black] (2,0) circle (0.5); 
				\node[circle,fill=black,draw=black,inner sep=1.5pt] (v1) at (-0.1,0) {};
				\node[circle,fill=black,draw=black,inner sep=1.5pt] (v2) at (1.42,0) {};
				\node at (0.65,0.3) {$P$};
				\draw [dash pattern=on 3pt off 2pt, line width=1.1pt](v1) -- (v2);
				\node at (0.65,-1.95) {\large$G_5^3$};
			\end{tikzpicture}
		};
		\node[right=0.5cm of fig5] (fig6) {
			\begin{tikzpicture}[scale =1]
				\draw[line width=1.1pt, thick, black] (-1.5,0) circle (0.5);
				\draw[line width=1.1pt, thick, black] (1.5,0) circle (0.5);
				\draw[line width=1.1pt, thick, black] (0,-1.5) circle (0.5); 
				\node[circle,fill=black,draw=black,inner sep=1.5pt] (v1) at (-0.1,0) {};
				\node[circle,fill=black,draw=black,inner sep=1.5pt] (v2) at (0.92,0) {};
				\node[circle,fill=black,draw=black,inner sep=1.5pt] (v3) at (-1.09,0) {};
				\node[circle,fill=black,draw=black,inner sep=1.5pt] (v4) at (-0.1,-1) {};
				\node at (0.22,0.3) {$P^1$};
				\node at (-0.77,0.3) {$P^2$};
				\node at (0.05,-0.52) {$P^3$};
				\draw [dash pattern=on 3pt off 2pt, line width=1.1pt](v1) -- (v2);
				\draw [dash pattern=on 3pt off 2pt, line width=1.1pt](v1) -- (v4);
				\draw [dash pattern=on 3pt off 2pt, line width=1.1pt](v1) -- (v3);
				\node at (-0.3,-2.5) {\large$G_6^3$};
				
			\end{tikzpicture}
		};
		\node[right=0.5cm of fig6] (fig7) {
			\begin{tikzpicture}[scale =1]
				\draw[line width=1.1pt, thick, black] (0,0) circle (0.5);
				\draw[line width=1.1pt, thick, black] (2,1.5) circle (0.5);
				\draw[line width=1.1pt, thick, black] (2,-1.5) circle (0.5); 
				\node[circle,fill=black,draw=black,inner sep=1.5pt] (v1) at (0.41,0) {};
				\node[circle,fill=black,draw=black,inner sep=1.5pt] (v2) at (1.42,1.5) {};
				\node[circle,fill=black,draw=black,inner sep=1.5pt] (v4) at (1.42,-1.5) {};
				\draw [dash pattern=on 3pt off 2pt, line width=1.1pt](0.5,0) -- (0.5,1.5);
				\draw [dash pattern=on 3pt off 2pt, line width=1.1pt](0.5,1.5) -- (1.45,1.5);
				\draw [dash pattern=on 3pt off 2pt, line width=1.1pt](0.5,0) -- (0.5,-1.5);
				\draw [dash pattern=on 3pt off 2pt, line width=1.1pt](0.5,-1.5) -- (1.45,-1.5);
				\node at (0.7,1.2) {$P^1$};
				\node at (0.7,-1.2) {$P^2$};
				\node at (0.8,-2.3) {\large$G_7^3$};

			\end{tikzpicture}
		};
      \node[below=0.8cm of fig4](fig8) {	  	
		\begin{tikzpicture}[scale =1]
			\hspace{-0.2cm}
			\draw[line width=1.3pt, thick, black] (0,0) circle (0.6);
			\draw[line width=1.3pt, thick, black] (-1.2,0) circle (0.6);
			\node[circle,fill=black,draw=black,inner sep=1.5pt] (v1) at (-0.6,0.05) {};
			\node[circle,fill=black,draw=black,inner sep=1.5pt] (v2) at (-0.9,0.6) {};
			\node[circle,fill=black,draw=black,inner sep=1.5pt] (v3) at (-1.8,0) {};
			\node at (-1.2,0.25) {$P$};
			\draw [dash pattern=on 3pt off 2pt, line width=1.5pt] (v2) -- (v3);
			\node at (-0.7,-1) {\large$G_1^4$};
		\end{tikzpicture}
	};
	\node[right=0.5cm of fig8] (fig9) {
		\begin{tikzpicture}[scale =1]
			\draw[line width=1.3pt, thick, black] (0,0) circle (0.6);
			\draw[line width=1.3pt, thick, black] (-1.2,0) circle (0.6);
			\node[circle,fill=black,draw=black,inner sep=1.5pt] (v1) at (-0.68,0) {};
			\node[circle,fill=black,draw=black,inner sep=1.5pt] (v3) at (-1.87,0) {};
			\draw [dash pattern=on 3pt off 2pt, line width=1.5pt] (v1) -- (v3);
			\node at (-1.45,-0.3) {$P$};
			\node at (-1,-1.4) {\large$G_2^4$};
		\end{tikzpicture}
	};
	\node[right=0.8cm of fig9] (fig10) {
		\begin{tikzpicture}[scale =1]
			\draw[line width=1.3pt, thick, black] (1,0) circle (0.6);
			\draw[line width=1.3pt, thick, black] (-1,0) circle (0.6);
			\node[circle,fill=black,draw=black,inner sep=1.5pt] (v1) at (0.33,0) {};
			\node[circle,fill=black,draw=black,inner sep=1.5pt] (v2) at (-0.48,0) {};
			\node[circle,fill=black,draw=black,inner sep=1.5pt] (v3) at (-1.05,-0.6) {};
			\node[circle,fill=black,draw=black,inner sep=1.5pt] (v4) at (-1.05,0.6) {};
			\draw [dash pattern=on 3pt off 2pt, line width=1.5pt] (v1) -- (v2);
			\draw [dash pattern=on 3pt off 2pt, line width=1.5pt] (v3) -- (v4);
			\node at (-0.25,0.25) {$P^2$};
			\node at (-1.05,0) {$P^1$};
			\node at (-0.25,-1.4) {\large$G_3^4$};
		\end{tikzpicture}
	};
	\node[right=0.8cm of fig10] (fig11) {
		\begin{tikzpicture}[scale =1]
			\draw[line width=1.3pt, thick, black] (1,0) circle (0.6);
			\draw[line width=1.3pt, thick, black] (-1,0) circle (0.6);
			\node[circle,fill=black,draw=black,inner sep=1.5pt] (v1) at (0.33,0) {};
			\node[circle,fill=black,draw=black,inner sep=1.5pt] (v2) at (-0.48,0) {};
			\node[circle,fill=black,draw=black,inner sep=1.5pt] (v3) at (-1.67,0) {};
			\draw [dash pattern=on 3pt off 2pt, line width=1.5pt] (v1) -- (v2);
			\draw [dash pattern=on 3pt off 2pt, line width=1.5pt] (v2) -- (v3);
			\node at (-0.25,0.25) {$P^2$};
			\node at (-1.25,0.25) {$P^1$};
			\node at (-0.25,-1.4) {\large$G_4^4$};
		\end{tikzpicture}
	};
\node [below=0.6cm of fig8](fig12) {	  	
	\begin{tikzpicture}[scale =1]
		\hspace{-0.2cm}
		\draw[line width=1.3pt, thick, black] (0,0) circle (0.8);
		\draw[line width=1.3pt, thick, black] (-0.8,0) circle (0.8);
		
		\node at (-0.5,-1.2) {\large$G_1^6$};
	\end{tikzpicture}
};
\node[right=0.9cm of fig12] (fig13) {
	\begin{tikzpicture}[scale =1]
		\draw[line width=1.3pt, thick, black] (0,0) circle (0.8);
		\node[circle,fill=black,draw=black,inner sep=1.5pt] (v1) at (-0.85,-0.2) {};
		\node[circle,fill=black,draw=black,inner sep=1.5pt] (v2) at (0.69,-0.2) {};
		\node[circle,fill=black,draw=black,inner sep=1.5pt] (v3) at (-0.1,0.8) {};
		\draw[dash pattern=on 3pt off 2pt, line width=1.1pt,black] (v1) to (v3);
		\draw[dash pattern=on 3pt off 2pt, line width=1.1pt,black] (v3) to (v2);
		\node at (-0.35,-1.42) {\large$G_2^6$};
	\end{tikzpicture}
};
\node[right=1.4cm of fig13] (fig14) {
	\begin{tikzpicture}[scale =1]
		\draw[line width=1.3pt, thick, black] (0,0) circle (0.8);
		\node[circle,fill=black,draw=black,inner sep=1.5pt] (v1) at (-0.76,0.4) {};
		\node[circle,fill=black,draw=black,inner sep=1.5pt] (v2) at (-0.77,-0.4) {};
		\node[circle,fill=black,draw=black,inner sep=1.5pt] (v3) at (0.61,0.4) {};
		\node[circle,fill=black,draw=black,inner sep=1.5pt] (v4) at (0.62,-0.4) {};
		\draw[dash pattern=on 3pt off 2pt, line width=1.1pt,black] (v1) to[bend right=90, looseness=2] (v2);
		\draw[dash pattern=on 3pt off 2pt, line width=1.1pt,black] (v3) to[bend left=90, looseness=2] (v4);
		\node at (-0.26,-1.47) {\large$G_3^6$};
	\end{tikzpicture}
};
\node[right=1.4cm of fig14] (fig15) {
	\begin{tikzpicture}[scale =1]
		\draw[line width=1.3pt, thick, black] (0,0) circle (0.8);
		\node[circle,fill=black,draw=black,inner sep=1.5pt] (v1) at (-0.76,0.4) {};
		\node[circle,fill=black,draw=black,inner sep=1.5pt] (v2) at (-0.1,-0.1) {};
		\node[circle,fill=black,draw=black,inner sep=1.5pt] (v3) at (0.61,0.4) {};
		\node[circle,fill=black,draw=black,inner sep=1.5pt] (v4) at (-0.1,-0.8) {};
		\draw [dash pattern=on 3pt off 2pt, line width=1.1pt] (v1) -- (v2);
		\draw [dash pattern=on 3pt off 2pt, line width=1.1pt] (v3) -- (v2);
		\draw [dash pattern=on 3pt off 2pt, line width=1.1pt] (v2) -- (v4);
		\node at (-0.36,-1.45) {\large$G_1^7$};
	\end{tikzpicture}
};
\end{tikzpicture}
\caption{The graphs $G_j^3$ ($j$=1,...,7), $G_j^4$ ($j$=1,...,4), $G_j^6$ ($j$=1,...,3)  and graph $G_1^7$ } 
\label{figure:G3} 
\end{figure}
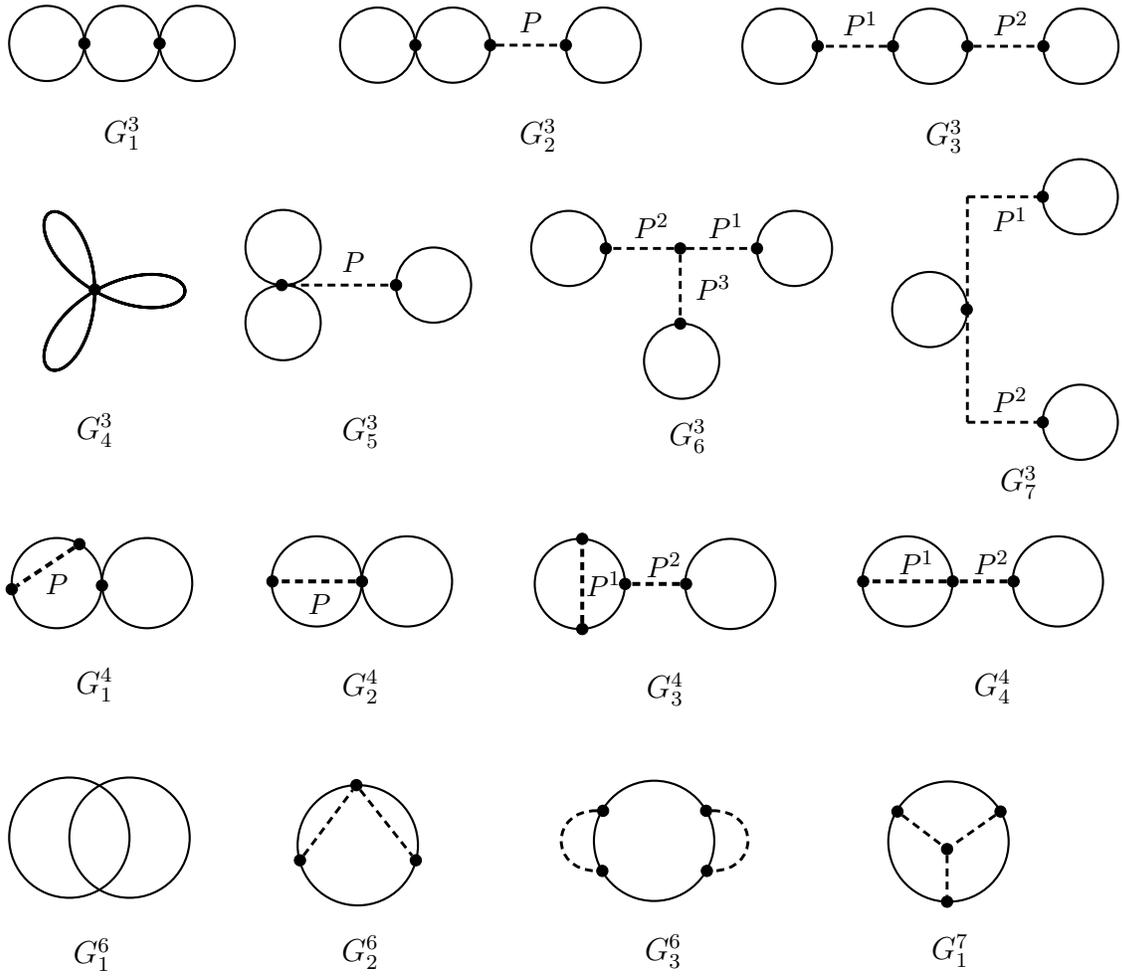

\section{Tricyclic graphs with \( \lambda_2(G) \) $<$ \(-\frac{1}{2}\) }
Let  $G$ is connected a tricyclic graph with the second largest distance eigenvalue less than \(-\frac{1}{2}\). According to Lemma \ref{th:ch-4}, $G$ is choral, then each cycle of $G$ is triangle. Similarly, each cycle of $B(G)$ is triangle. Next, by Lemma \ref{th:ch-5}, we will consider $B(G)$ in four parts.

When $B(G)\cong G_j^3$ ($j$=1,...,7), according to the symmetry of the triangle, $G_1^3$ and $G_2^3$ can be obtained through $G_3^3$. Similarly, $G_4^3$ and $G_5^3$ can be obtained through $G_7^3$. But if one of the three paths in $G_6^3$ has a length of zero, then $G_6^3$ becomes $G_7^3$. If the length of the three paths in $G_6^3$ is not zero, since $F_5$ is forbidden subgraph, then we can get that $G_6^3$ of this structure are not subgraphs in a tricyclic graph with the second largest distance eigenvalue less than \(-\frac{1}{2}\). 
When $B(G)\cong G_j^4$ ($j$=1,...,4),we can get that $G_1^4$, $G_3^4$ and $G_4^4$ are not subgraphs in a tricyclic graph we are considering, since $F_9$, $F_{10}$ are forbidden subgraphs. 
When $B(G)\cong G_j^6$ ($j$=1,2,3), we can get that $G_1^6$, $G_2^6$ and $G_3^6$ are not subgraphs in a tricyclic graph we are considering, since $F_{12}$ and $F_{11}$ are forbidden subgraphs and each cycle of $B(G)$ is triangle. When $B(G)\cong G_1^7$, then $B(G)\cong G_7^1$ is not considered in this article, since $K_4 \cong G_1^7$ and Guo and Zhou had characterized split graphs with $\lambda _2\left( G \right) \le -\frac{1}{2}$ in \cite{RF6}. Therefore, we consider $B(G)\cong G_3^3,G_7^3,G_2^4$ in the subsequent part and we give the following definitions.

A path \( u_0u_1\ldots u_r \) with \( r \ge 1 \) in a graph \( G \) is called a pendant path of length $r$ at $u_0$ if $d_G(u_0)\ge3$, the degrees of $u_1,u_2, \ldots, u_{r-1}$ (if any presents) are all equal to 2 in \( G \), and $d_G(u_r)=1$. In such a case, it is also said that $G$ is constructed from $G-\{u_1,u_2,\ldots,u_r\}$ by attaching a pendant path of length $r$ at \(u_0\). When attaching a pendant path of length 0 to a vertex $v$ in graph $G$, the graph remains unchanged. Additionally, a pendant path of length 1 at \( u_0 \) is called a pendant edge at \(u_0\).
 
 Let $s$ and $t$ be integers with $s\ge0$ and $t\ge0$. Denote by $T_1(s,t)$ the tricyclic graph obtained from the three triangles $\varDelta _{11}=u_1u_2u_3$, $\varDelta _{12}=v_1v_2v_3$ and $\varDelta _{13}=w_1w_2w_3$ and the two paths $P_{s}= p_0 p_1...p_s$ and $P_{t}=q_0q_1...q_t$ by identifying $u_3$ with $p_0$, identifying  $v_1$ with $p_s$, identifying $v_3$ with $q_0$ and identifying $w_1$ with $q_t$. The graph $T_1(s,t)$ is illustrated in Figure.\ref{figure:t1}.
 
\begin{figure}[htbp]
\centering
\begin{tikzpicture}[scale =2.5]
	\node[circle,fill=black,draw=black,inner sep=1.5pt] (v1) at (-0.3,0) {};
	\node[circle,fill=black,draw=black,inner sep=1.5pt] (v2) at (0.3,0) {};
	\node[circle,fill=black,draw=black,inner sep=1.5pt] (v3) at (0,0) {};
	\node[circle,fill=black,draw=black,inner sep=1.5pt] (v4) at (-0.6,0) {};
	\node[circle,fill=black,draw=black,inner sep=1.5pt] (v5) at (0.6,0) {};
	\node[circle,fill=black,draw=black,inner sep=1.5pt] (v6) at (-1.3,0) {};
	\node[circle,fill=black,draw=black,inner sep=1.5pt] (v7) at (1.3,0) {};
	\node[circle,fill=black,draw=black,inner sep=1.5pt] (v8) at (-1.6,0) {};
	\node[circle,fill=black,draw=black,inner sep=1.5pt] (v9) at (1.6,0) {};
	\node[circle,fill=black,draw=black,inner sep=1.5pt] (v10) at (-2.2,0) {};
	\node[circle,fill=black,draw=black,inner sep=1.5pt] (v11) at (2.2,0) {};
	\node[circle,fill=black,draw=black,inner sep=1.5pt] (v12) at (-1.9,0.6) {};
	\node[circle,fill=black,draw=black,inner sep=1.5pt] (v13) at (1.9,0.6) {};
	\node[circle,fill=black,draw=black,inner sep=1.5pt] (v14) at (0.3,0.6) {};
	\node[circle,fill=black,draw=black,inner sep=1.5pt] (v15) at (-0.3,0.6) {};
	\node[below] at (v2) {$v_{q-1}$};
	\node[below] at (v7) {$v_1$};
	\node[below] at (v9) {$v_0$};
	\node[below] at (v1) {$u_{p-1}$};
	\node[below] at (v3) {$y_1$};
	\node[below] at (v6) {$u_1$};
	\node[below] at (v8) {$u_0$};
	
	\fill (-0.85,0) circle (0.8pt); 
	\fill (-0.95,0) circle (0.8pt); 
	\fill (-1.05,0) circle (0.8pt); 
	\fill (0.85,0) circle (0.8pt); 
	\fill (0.95,0) circle (0.8pt); 
	\fill (1.05,0) circle (0.8pt); 
	
	\draw [line width=1pt](0.6,0) -- (0.75,0);
	\draw[line  width=1pt] (1.15,0) -- (1.3,0);
	\draw [line width=1pt](-0.6,0) -- (-0.75,0);
	\draw[line  width=1pt] (-1.15,0) -- (-1.3,0);

	\draw [line width=1pt](v1) -- (v2);
	\draw[line  width=1pt] (v1) -- (v3);
	\draw [line width=1pt](v1) -- (v4);
	\draw [line width=1pt](v2) -- (v3);
	\draw [line width=1pt](v2) -- (v5);
	\draw [line width=1pt](v7) -- (v9);
	\draw [line width=1pt](v9) -- (v11);
	\draw [line width=1pt](v9) -- (v13);
	\draw [line width=1pt](v11) -- (v13);
	\draw [line width=1pt](v6) -- (v8);
	\draw [line width=1pt](v8) -- (v10);
	\draw [line width=1pt](v8) -- (v12);
	\draw [line width=1pt](v10) -- (v12);
	\draw [line width=1pt](v3) -- (v14);
	\draw [line width=1pt](v3) -- (v15);
	\draw [line width=1pt](v14) -- (v15);				
\end{tikzpicture}
\caption{Garph $T_2(p,q)$ } 
\label{figure:t2} 
\end{figure}
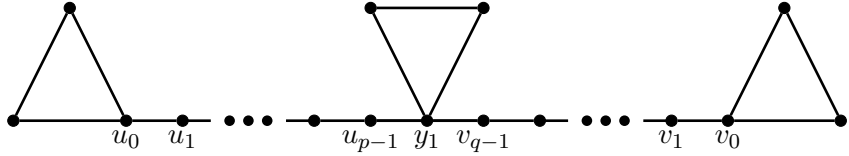

 Let $p$ and $q$ be integers with $p\ge 0$ and $q\ge 0$. Denote by $T_2(p,q)$ the tricyclic graph obtained from the three triangles $\varDelta _{21}=x_1x_2x_3$, $\varDelta _{22}=y_1y_2y_3$ and $\varDelta _{23}=z_1z_2z_3$ and the two paths $P_{p}= u_0 u_1...u_p$ and $P_{q}=v_0v_1...v_q$ by identifying $x_3$ with $u_0$, identifying $y_1$ with $u_p$, identifying  $y_1$ with $v_q$ and identifying  $z_1$ with $v_0$. The graph $T_2(p,q)$ is illustrated in Figure.\ref{figure:t2}.
  
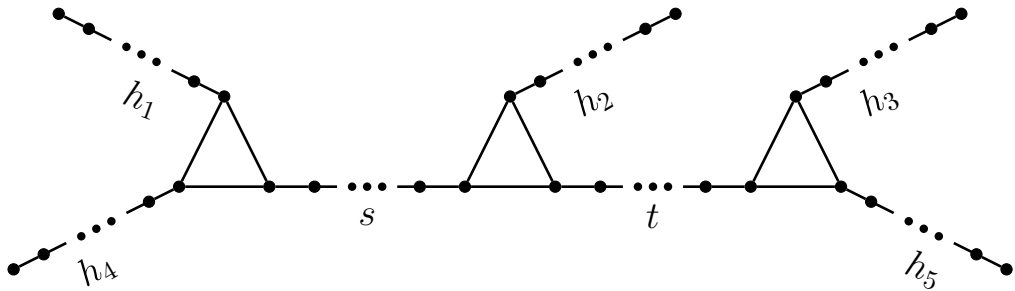
\begin{figure}[htbp]
	\centering
	\begin{tikzpicture}[scale =2]
		\node[circle,fill=black,draw=black,inner sep=1.5pt] (v1) at (-0.3,0) {};
		\node[circle,fill=black,draw=black,inner sep=1.5pt] (v2) at (0.3,0) {};
		\node[circle,fill=black,draw=black,inner sep=1.5pt] (v3) at (0,0.6) {};
		\node[circle,fill=black,draw=black,inner sep=1.5pt] (v4) at (-0.6,0) {};
		\node[circle,fill=black,draw=black,inner sep=1.5pt] (v5) at (0.6,0) {};
		\node[circle,fill=black,draw=black,inner sep=1.5pt] (v6) at (-1.3,0) {};
		\node[circle,fill=black,draw=black,inner sep=1.5pt] (v7) at (1.3,0) {};
		\node[circle,fill=black,draw=black,inner sep=1.5pt] (v8) at (-1.6,0) {};
		\node[circle,fill=black,draw=black,inner sep=1.5pt] (v9) at (1.6,0) {};
		\node[circle,fill=black,draw=black,inner sep=1.5pt] (v10) at (-2.2,0) {};
		\node[circle,fill=black,draw=black,inner sep=1.5pt] (v11) at (2.2,0) {};
		\node[circle,fill=black,draw=black,inner sep=1.5pt] (v12) at (-1.9,0.6) {};
		\node[circle,fill=black,draw=black,inner sep=1.5pt] (v13) at (1.9,0.6) {};
		\node[circle,fill=black,draw=black,inner sep=1.5pt] (v14) at (0.2,0.7) {};
		\node[circle,fill=black,draw=black,inner sep=1.5pt] (v15) at (0.9,1.05) {};
		\node[circle,fill=black,draw=black,inner sep=1.5pt] (v16) at (1.1,1.15) {};
		\node[circle,fill=black,draw=black,inner sep=1.5pt] (v17) at (2.1,0.7) {};
		\node[circle,fill=black,draw=black,inner sep=1.5pt] (v18) at (2.8,1.05) {};
		\node[circle,fill=black,draw=black,inner sep=1.5pt] (v19) at (3.0,1.15) {};
		\node[circle,fill=black,draw=black,inner sep=1.5pt] (v20) at (-2.1,0.7) {};
		\node[circle,fill=black,draw=black,inner sep=1.5pt] (v21) at (-2.8,1.05) {};
		\node[circle,fill=black,draw=black,inner sep=1.5pt] (v22) at (-3.0,1.15) {};
		\node[circle,fill=black,draw=black,inner sep=1.5pt] (v23) at (-2.4,-0.1) {};
		\node[circle,fill=black,draw=black,inner sep=1.5pt] (v24) at (-3.1,-0.45) {};
		\node[circle,fill=black,draw=black,inner sep=1.5pt] (v25) at (-3.3,-0.55) {};
		\node[circle,fill=black,draw=black,inner sep=1.5pt] (v26) at (2.4,-0.1) {};
		\node[circle,fill=black,draw=black,inner sep=1.5pt] (v27) at (3.1,-0.45) {};
		\node[circle,fill=black,draw=black,inner sep=1.5pt] (v28) at (3.3,-0.55) {};
		\node [rotate=-20] at (-2.45,0.58) {\Large$h_1$};
		\node [rotate=30] at (0.55,0.58) {\Large$h_2$};
		\node [rotate=30] at (2.45,0.58) {\Large$h_3$};
		\node [rotate=30] at (-2.75,-0.55) {\Large$h_4$};
		\node [rotate=-20] at (2.75,-0.55) {\Large$h_5$};
		\node at (-0.95,-0.2) {\Large$s$};
		\node at (0.95,-0.2) {\Large$t$};
		
		\fill (-0.85,0) circle (0.8pt); 
		\fill (-0.95,0) circle (0.8pt); 
		\fill (-1.05,0) circle (0.8pt); 
		\fill (0.85,0) circle (0.8pt); 
		\fill (0.95,0) circle (0.8pt); 
		\fill (1.05,0) circle (0.8pt); 
		
		\fill (0.45,0.83) circle (0.8pt); 
		\fill (0.55,0.88) circle (0.8pt); 
		\fill (0.65,0.93) circle (0.8pt); 
		
		\draw [line width=1pt](0.2,0.7) -- (0.35,0.775);
		\draw[line  width=1pt] (0.75,0.975) -- (0.9,1.05);
		
		\draw [line width=1pt](2.1,0.7) -- (2.25,0.775);
		\draw[line  width=1pt] (2.65,0.975) -- (2.8,1.05);
		
		\fill (2.35,0.83) circle (0.8pt); 
		\fill (2.45,0.88) circle (0.8pt); 
		\fill (2.55,0.93) circle (0.8pt); 
		
		\draw [line width=1pt](-2.1,0.7) -- (-2.25,0.775);
		\draw[line  width=1pt] (-2.65,0.975) -- (-2.8,1.05);
		
		\fill (-2.35,0.83) circle (0.8pt); 
		\fill (-2.45,0.88) circle (0.8pt); 
		\fill (-2.55,0.93) circle (0.8pt); 
		
		\draw [line width=1pt](-2.4,-0.1) -- (-2.55,-0.175);
		\draw[line  width=1pt] (-2.95,-0.375) -- (-3.1,-0.45);
		
		\fill (-2.65,-0.23) circle (0.8pt); 
		\fill (-2.75,-0.28) circle (0.8pt); 
		\fill (-2.85,-0.33) circle (0.8pt); 
		
		\draw [line width=1pt](2.4,-0.1) -- (2.55,-0.175);
		\draw[line  width=1pt] (2.95,-0.375) -- (3.1,-0.45);
		
		\fill (2.65,-0.23) circle (0.8pt); 
		\fill (2.75,-0.28) circle (0.8pt); 
		\fill (2.85,-0.33) circle (0.8pt); 

		\draw [line width=1pt](0.6,0) -- (0.75,0);
		\draw[line  width=1pt] (1.15,0) -- (1.3,0);
		\draw [line width=1pt](-0.6,0) -- (-0.75,0);
		\draw[line  width=1pt] (-1.15,0) -- (-1.3,0);

		\draw [line width=1pt](v1) -- (v2);
		\draw[line  width=1pt] (v1) -- (v3);
		\draw [line width=1pt](v1) -- (v4);
		\draw [line width=1pt](v2) -- (v3);
		\draw [line width=1pt](v2) -- (v5);
		\draw [line width=1pt](v7) -- (v9);
		\draw [line width=1pt](v9) -- (v11);
		\draw [line width=1pt](v9) -- (v13);
		\draw [line width=1pt](v11) -- (v13);
		\draw [line width=1pt](v6) -- (v8);
		\draw [line width=1pt](v8) -- (v10);
		\draw [line width=1pt](v8) -- (v12);
		\draw [line width=1pt](v10) -- (v12);
		\draw [line width=1pt](v3) -- (v14);
		\draw [line width=1pt](v15) -- (v16);
		\draw [line width=1pt](v13) -- (v17);
		\draw [line width=1pt](v18) -- (v19);
		\draw [line width=1pt](v12) -- (v20);
		\draw [line width=1pt](v21) -- (v22);	
		\draw [line width=1pt](v10) -- (v23);
		\draw [line width=1pt](v24) -- (v25);	
		\draw [line width=1pt](v11) -- (v26);
		\draw [line width=1pt](v27) -- (v28);		
	\end{tikzpicture}
	\caption{Graphs $T(s,t;h_1,h_2,h_3,h_4,h_5)$ } 
	\label{figure:tt} 
\end{figure}

  Let $u_1,u_2,u_3,u_4,u_5$ be the five vertices with degree two in the three triangles of  $T_1(s,t)$ with $u_1u_2$, $u_4u_5\in E(T_1(s,t))$, where $s\ge 0$ and $t\ge 0$. We use $T(s,t;h_1,h_2,h_3,h_4,h_5)$ to denote the graph obtained from $T_1(s,t)$ by attaching a pendant path of length $h_i$ at $u_i$, respectively, where $h_i$ $\ge 0$ for $1\le i\le 5$, as illustrated in Figure.\ref{figure:tt}.
  
\begin{figure}[h]
	\centering
	\begin{tikzpicture}
		\node (fig1) {	  	
			\begin{tikzpicture}[scale =2]
				\node[circle,fill=black,draw=black,inner sep=2pt] (v1) at (0,0) {};
				\node[circle,fill=black,draw=black,inner sep=2pt] (v2) at (0.9,0) {};
				\node[circle,fill=black,draw=black,inner sep=2pt] (v3) at (0.45,0.8) {};
				\node[circle,fill=black,draw=black,inner sep=2pt] (v4) at (-0.9,0) {};
				\node[circle,fill=black,draw=black,inner sep=2pt] (v5) at (-0.45,0.8) {};
				\node[circle,fill=black,draw=black,inner sep=2pt] (v6) at (-0.7,-0.6) {};
				\node[circle,fill=black,draw=black,inner sep=2pt] (v7) at (0.15,-0.85) {};
				\node[circle,fill=black,draw=black,inner sep=2pt] (v8) at (0.4,-0.65) {};
				\node[circle,fill=black,draw=black,inner sep=2pt] (v9) at (0.75,-0.2) {};
				\node at (0.1,-1.2) {\Large$T_3^{k}$};
				\node[rotate=-30] at (0.85,-0.65) {\Large$k$};
				\fill (0.5,-0.55) circle (1.1pt); 
				\fill (0.59,-0.43) circle (1.1pt); 
				\fill (0.68,-0.32) circle (1.1pt); 
				\draw [decorate, decoration={brace, amplitude=7pt,raise=6pt,mirror},black,very thick](0.4,-0.65) -- (0.75,-0.2) ;
				\draw [line width=1.5pt](v1) -- (v2);
				\draw[line width=1.5pt] (v1) -- (v3);
				\draw [line width=1.5pt](v1) -- (v4);
				\draw [line width=1.5pt](v1) -- (v5);
				\draw [line width=1.5pt](v1) -- (v6);
				\draw [line width=1.5pt](v1) -- (v7);
				\draw [line width=1.5pt](v2) -- (v3);
				\draw [line width=1.5pt](v4) -- (v5);
				\draw [line width=1.5pt](v1) -- (v8);
				\draw [line width=1.5pt](v1) -- (v9);
				\draw [line width=1.5pt](v6) -- (v7);
				\draw [line width=1.5pt](v6) -- (v7);
				
			\end{tikzpicture}
		};
		\node[right=3.5cm of fig1] (fig2) {
			\begin{tikzpicture}[scale =2]
				\node[circle,fill=black,draw=black,inner sep=2pt] (v1) at (0,0) {};
				\node[circle,fill=black,draw=black,inner sep=2pt] (v2) at (0.7,0.7) {};
				\node[circle,fill=black,draw=black,inner sep=2pt] (v3) at (0.7,-0.7) {};
				\node[circle,fill=black,draw=black,inner sep=2pt] (v4) at (-0.5,0.7) {};
				\node[circle,fill=black,draw=black,inner sep=2pt] (v5) at (-1,0) {};
				\node[circle,fill=black,draw=black,inner sep=2pt] (v6) at (-0.5,-0.7) {};
				\node[circle,fill=black,draw=black,inner sep=2pt] (v7) at (-0.2,-0.7) {};
				\node[circle,fill=black,draw=black,inner sep=2pt] (v8) at (0.4,-0.7) {};
				\node at (0.1,-1.35) {\Large$T_4^{t}$};
				\node[rotate=0] at (0.05,-1) {\Large$t$};
				\fill (0.0,-0.7) circle (1.1pt); 
				\fill (0.15,-0.7) circle (1.1pt); 
				\fill (0.3,-0.7) circle (1.1pt); 
				\draw [decorate, decoration={brace, amplitude=6pt,raise=5pt,mirror},black,very thick](-0.15,-0.7) -- (0.45,-0.7) ;
				\draw [line width=1.5pt](v1) -- (v2);
				\draw[line width=1.5pt] (v1) -- (v3);
				\draw [line width=1.5pt](v1) -- (v4);
				\draw [line width=1.5pt](v1) -- (v5);
				\draw [line width=1.5pt](v1) -- (v6);
				\draw [line width=1.5pt](v1) -- (v7);
				\draw [line width=1.5pt](v2) -- (v3);
				\draw [line width=1.5pt](v5) -- (v6);
				\draw [line width=1.5pt](v1) -- (v8);
				\draw [line width=1.5pt](v4) -- (v5);
			\end{tikzpicture}
		};
	\end{tikzpicture}
	\caption{Graphs $T_{3}^{k}$ and $T_{4}^{t}$ } 
	\label{fig:5} 
\end{figure}
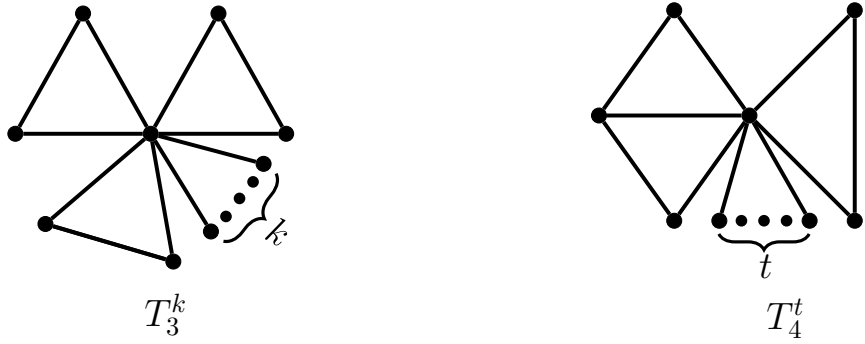

  Let $k$ be a interger with $k\ge 0$. Denote by $T_{3}^{k}$ the tricyclic graph obtained from the $T_2(0,0)$ and attaching $k$ pendant edges at the vertex with maximum degree. The graph $T_{3}^{k}$ is illustrated in Figure.\ref{fig:5}.
  
   Let $t$ be a interger with $t\ge 0$. Denote by $T_{4}^{t}$ the tricyclic graph obtained from the three triangles $\varDelta_{41}=x_1x_2x_3$, $\varDelta_{42}=y_1y_2y_3$ and $\varDelta_{43}=z_1z_2z_3$ by identifying $x_2$ with $y_2$, identifying $x_3$ with $y_3$, identifying $y_2$ with $z_1$ and attaching $t$ pendant edges at the vertex with maximum degree. The graph $T_{4}^{t}$ is illustrated in  Figure.\ref{fig:5}.
\begin{figure}[htbp]
	\centering
	\begin{tikzpicture}[scale =2.5]
		\node[circle,fill=black,draw=black,inner sep=1.5pt] (v1) at (-0.3,0) {};
		\node[circle,fill=black,draw=black,inner sep=1.5pt] (v2) at (0.3,0) {};
		\node[circle,fill=black,draw=black,inner sep=1.5pt] (v3) at (0,0.6) {};
		\node[circle,fill=black,draw=black,inner sep=1.5pt] (v4) at (-0.6,0) {};
		\node[circle,fill=black,draw=black,inner sep=1.5pt] (v5) at (0.6,0) {};
		\node[circle,fill=black,draw=black,inner sep=1.5pt] (v6) at (-1.3,0) {};
		\node[circle,fill=black,draw=black,inner sep=1.5pt] (v7) at (1.3,0) {};
		\node[circle,fill=black,draw=black,inner sep=1.5pt] (v8) at (-1.6,0) {};
		\node[circle,fill=black,draw=black,inner sep=1.5pt] (v9) at (1.6,0) {};
		\node[circle,fill=black,draw=black,inner sep=1.5pt] (v10) at (-2.2,0) {};
		\node[circle,fill=black,draw=black,inner sep=1.5pt] (v11) at (2.2,0) {};
		\node[circle,fill=black,draw=black,inner sep=1.5pt] (v12) at (-1.9,0.6) {};
		\node[circle,fill=black,draw=black,inner sep=1.5pt] (v13) at (1.9,0.6) {};
		\node (v14) at (0,-0.15) {};
		\node[below] at (v2) {$q_0$};
		\node[below] at (v5) {$q_1$};
		\node[below] at (v7) {$q_{t-1}$};
		\node[below] at (v9) {$q_t$};
		\node[below] at (v1) {$p_0$};
		\node[below] at (v4) {$p_1$};
		\node[below] at (v6) {$p_{s-1}$};
		\node[below] at (v8) {$p_s$};
		
		\fill (-0.85,0) circle (0.8pt); 
		\fill (-0.95,0) circle (0.8pt); 
		\fill (-1.05,0) circle (0.8pt); 
		\fill (0.85,0) circle (0.8pt); 
		\fill (0.95,0) circle (0.8pt); 
		\fill (1.05,0) circle (0.8pt); 
		
		\draw [line width=1pt](0.6,0) -- (0.75,0);
		\draw[line  width=1pt] (1.15,0) -- (1.3,0);
		\draw [line width=1pt](-0.6,0) -- (-0.75,0);
		\draw[line  width=1pt] (-1.15,0) -- (-1.3,0);

		\draw [line width=1pt](v1) -- (v2);
		\draw[line  width=1pt] (v1) -- (v3);
		\draw [line width=1pt](v1) -- (v4);
		\draw [line width=1pt](v2) -- (v3);
		\draw [line width=1pt](v2) -- (v5);
		\draw [line width=1pt](v7) -- (v9);
		\draw [line width=1pt](v9) -- (v11);
		\draw [line width=1pt](v9) -- (v13);
		\draw [line width=1pt](v11) -- (v13);
		\draw [line width=1pt](v6) -- (v8);
		\draw [line width=1pt](v8) -- (v10);
		\draw [line width=1pt](v8) -- (v12);
		\draw [line width=1pt](v10) -- (v12);			
	\end{tikzpicture}
	\caption{Garph $T_1(s,t)$ } 
	\label{figure:t1} 
\end{figure}
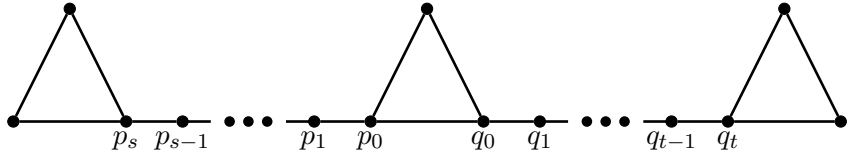

\noindent\begin{theorem}\label{th:ch-6}
	Let $G$ be a tricyclic graph. Then $\lambda_2(G)<-\frac{1}{2}$ if and only if
	\begin{itemize}
		\item $G \cong T_{5},T_{6},T_{7}$ as displayed in Figure.\ref{figure:t5}
		\item  $G\cong T(s,t;h_1,h_2,h_3,h_4,h_5)$, where $s\ge 0$, $t\ge 0$ and $h_i\ge 0$ for $1\le i\le 5$, or
		\item $G\cong T_{3}^{k}$ with $k\ge 0$, or
		\item $G\cong T_{4}^{t}$, where $t\ge 0$.
	\end{itemize}
\end{theorem}
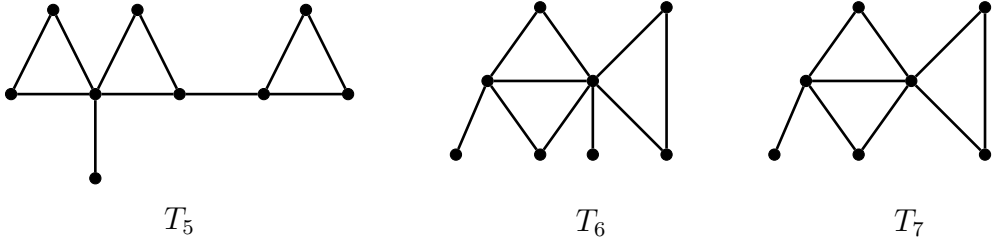
\begin{figure}[h]
	\centering
		\centering
	\begin{tikzpicture}
		\node (fig1) {	  	
			\begin{tikzpicture}[scale =1.4]
	
				\node[circle,fill=black,draw=black,inner sep=1.5pt] (v1) at (0,0) {};
				\node[circle,fill=black,draw=black,inner sep=1.5pt] (v2) at (0.8,0) {};
				\node[circle,fill=black,draw=black,inner sep=1.5pt] (v3) at (1.6,0) {};
				\node[circle,fill=black,draw=black,inner sep=1.5pt] (v4) at (1.2,0.8) {};
				\node[circle,fill=black,draw=black,inner sep=1.5pt] (v5) at (-0.8,0) {};
				\node[circle,fill=black,draw=black,inner sep=1.5pt] (v6) at (-1.6,0) {};
				\node[circle,fill=black,draw=black,inner sep=1.5pt] (v7) at (-1.2,0.8) {};
				\node[circle,fill=black,draw=black,inner sep=1.5pt] (v8) at (-0.4,0.8) {};
				\node[circle,fill=black,draw=black,inner sep=1.5pt] (v9) at (-0.8,-0.8){};
				\node at (0.0,-1.2) {\large$T_5$};
				
				\draw [line width=1pt](v1) -- (v2);
				\draw[line width=1pt] (v1) -- (v5);
				\draw [line width=1pt](v2) -- (v4);
				\draw [line width=1pt](v2) -- (v3);
				\draw [line width=1pt](v3) -- (v4);
				\draw [line width=1pt](v1) -- (v8);
				\draw [line width=1pt](v5) -- (v8);
				\draw [line width=1pt](v5) -- (v9);
				\draw [line width=1pt](v5) -- (v6);
				\draw [line width=1pt](v5) -- (v7);
				\draw [line width=1pt](v6) -- (v7);
				
			\end{tikzpicture}
		};
		\node[right=1cm of fig1] (fig2) {
			\begin{tikzpicture}[scale =1.4]
	
				\node[circle,fill=black,draw=black,inner sep=1.5pt] (v1) at (0,0) {};
				\node[circle,fill=black,draw=black,inner sep=1.5pt] (v2) at (0.7,0.7) {};
				\node[circle,fill=black,draw=black,inner sep=1.5pt] (v3) at (0.7,-0.7) {};
				\node[circle,fill=black,draw=black,inner sep=1.5pt] (v4) at (-0.5,0.7) {};
				\node[circle,fill=black,draw=black,inner sep=1.5pt] (v5) at (-1,0) {};
				\node[circle,fill=black,draw=black,inner sep=1.5pt] (v6) at (-0.5,-0.7) {};
				\node[circle,fill=black,draw=black,inner sep=1.5pt] (v7) at (0,-0.7) {};
				\node[circle,fill=black,draw=black,inner sep=1.5pt] (v8) at (-1.3,-0.7) {};
				\node at (-0.2,-1.35) {\large$T_6$};
				
				\draw [line width=1pt](v1) -- (v2);
				\draw[line width=1pt] (v1) -- (v3);
				\draw [line width=1pt](v1) -- (v4);
				\draw [line width=1pt](v1) -- (v5);
				\draw [line width=1pt](v1) -- (v6);
				\draw [line width=1pt](v1) -- (v7);
				\draw [line width=1pt](v2) -- (v3);
				\draw [line width=1pt](v5) -- (v6);
				\draw [line width=1pt](v5) -- (v8);
				\draw [line width=1pt](v4) -- (v5);
			\end{tikzpicture}
		};
		\node[right=1cm of fig2] (fig3) {
			\begin{tikzpicture}[scale =1.4]
				\node[circle,fill=black,draw=black,inner sep=1.5pt] (v1) at (0,0) {};
				\node[circle,fill=black,draw=black,inner sep=1.5pt] (v2) at (0.7,0.7) {};
				\node[circle,fill=black,draw=black,inner sep=1.5pt] (v3) at (0.7,-0.7) {};
				\node[circle,fill=black,draw=black,inner sep=1.5pt] (v4) at (-0.5,0.7) {};
				\node[circle,fill=black,draw=black,inner sep=1.5pt] (v5) at (-1,0) {};
				\node[circle,fill=black,draw=black,inner sep=1.5pt] (v6) at (-0.5,-0.7) {};
				\node[circle,fill=black,draw=black,inner sep=1.5pt] (v8) at (-1.3,-0.7) {};
				\node at (-0.2,-1.35) {\large$T_7$};
				
				\draw [line width=1pt](v1) -- (v2);
				\draw[line width=1pt] (v1) -- (v3);
				\draw [line width=1pt](v1) -- (v4);
				\draw [line width=1pt](v1) -- (v5);
				\draw [line width=1pt](v1) -- (v6);
				\draw [line width=1pt](v2) -- (v3);
				\draw [line width=1pt](v5) -- (v6);
				\draw [line width=1pt](v5) -- (v8);
				\draw [line width=1pt](v4) -- (v5);
			\end{tikzpicture}
		};
	\end{tikzpicture}
	\caption{Graphs $ T_{5} ,T_{6}$ and $ T_{7}$ }
	\label{figure:t5} 
\end{figure}

\begin{Tproof}\textbf{~of~Theorem~\ref{th:ch-6}.}
Let  $G$ is connected  tricyclic graph for which the second largest distance eigenvalue is less than \(-\frac{1}{2}\). By Lemma \ref{th:ch-4}, $G$ is choral, then each cycle of $G$ is triangle. Next, we consider the following three cases based on whether the  tricyclic graphs $G$ are based on $T_1(s,t)$ graphs or $T_2(p,q)$ graphs or $T_{4}^{t}$ graphs.

\noindent\textbf{Case 1.} When $G$ are based on $T_1(s,t)$ graph, we discuss the following three subcases.\\
\textbf{Subcase 1.1.} $s=0$ and $t=0$.\

In this case, since $F_{5}$ and $F_{8}$ are forbidden subgraphs, then a pendant path may exist at the vertices of degree two in the three triangles of $T_1(0,0)$. Given that $F_{5}$ is forbidden subgraph, the degree of the vertex which is not on the three triangles is at most two. Consequently, $G\cong T(0,0;h_1,h_2,h_3,h_4,h_5)$, where $h_i\ge 0$ for $1\le i\le 5$.

\begin{figure}[h]
	\centering
	\begin{tikzpicture}[scale =2]
		\node[circle,fill=black,draw=black,inner sep=2pt] (v1) at (0,0) {};
		\node[circle,fill=black,draw=black,inner sep=2pt] (v2) at (0.8,0) {};
		\node[circle,fill=black,draw=black,inner sep=2pt] (v3) at (1.6,0) {};
		\node[circle,fill=black,draw=black,inner sep=2pt] (v4) at (1.2,0.8) {};
		\node[circle,fill=black,draw=black,inner sep=2pt] (v5) at (-0.8,0) {};
		\node[circle,fill=black,draw=black,inner sep=2pt] (v6) at (-1.6,0) {};
		\node[circle,fill=black,draw=black,inner sep=2pt] (v7) at (-1.2,0.8) {};
		\node[circle,fill=black,draw=black,inner sep=2pt] (v8) at (-0.4,0.8) {};
		
		\node at (0.8,-0.2) {\Large$w_1$};
		\node at (1.6,-0.2) {\Large$w_3$};
		\node at (0.0,-0.2) {\Large$v_1$};
		\node at (-0.8,-0.2) {\Large$u_3(v_3)$};
		\node at (-1.6,-0.2) {\Large$u_1$};
		\node [right] at (1.2,0.8) {\Large$w_2$};
		\node [right] at (-1.2,0.8) {\Large$v_2$};
		\node[right] at (-0.4,0.8) {\Large$u_2$};
		\draw [line width=1.5pt](v1) -- (v2);
		\draw[line width=1.5pt] (v1) -- (v5);
		\draw [line width=1.5pt](v2) -- (v4);
		\draw [line width=1.5pt](v2) -- (v3);
		\draw [line width=1.5pt](v3) -- (v4);
		\draw [line width=1.5pt](v1) -- (v8);
		\draw [line width=1.5pt](v5) -- (v8);
		\draw [line width=1.5pt](v5) -- (v6);
		\draw [line width=1.5pt](v5) -- (v7);
		\draw [line width=1.5pt](v6) -- (v7);	
	\end{tikzpicture}
	\caption{Graph $T_1(1,0)$}
	\label{fig:7} 
\end{figure}
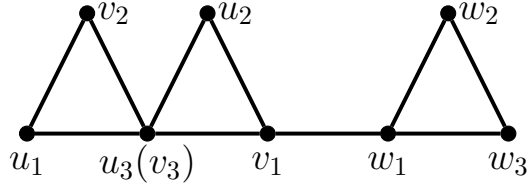
\noindent\textbf{Subcase 1.2.} $s=1$ or $t=1$.

For the sake of description, we label the vertices of $T_1(1,0)$ as shown in Figure.\ref{fig:7}. Then $T_1(1,0)$ are tricyclic garph obtained from the three 
triangles $\varDelta _{11}=u_1u_2u_3$, $\varDelta _{12}=v_1v_2v_3$ and $\varDelta _{13}=w_1w_2w_3$ by identifying $u_{3}$ with $v_{3}$ and  $v_{1}w_{1}\in E(T_1(1,0))$. 
Considering that $F_{5}$ is a forbidden subgraph, we obtain $d(v_1) =d(w_1) =3$. Due to that $F_1$ is a forbidden subgraph, we obtain $d(u_3)=4$ or 5. 

 Assume that  $d\left( u_3 \right) =5$. Let $N_G\left( u_3 \right)\setminus \{u_1,u_2,v_1,v_2\}=x_{1}$. We obtain $d(x_1)=1$, since $F_{7}$ is a forbidden subgraph. Because $F_{3}$ and $F_{4}$ are forbidden subgraphs, we obtain $d(u_1)=d(u_2)=d(w_2)=d(w_3)=2$. We obtain $d(v_2)=2$, since $F_{6}$ is a forbidden subgraph. Therefore $G\cong T_5$.
 
 Assume next that $d(u_3) =4$, since $F_6$ is a forbidden subgraph, we obtain $d(v_2)\le 3$. Because $F_4$ is a forbidden subgraph, we obtain $d(u_1),d(u_2),d(w_2),d(w_3)\le3$.  We obtain that there is a pendant path at the vertex of degree two in the three triangles of $T_1(1,0)$, since $F_{4}$ and $F_{5}$ are forbidden subgraphs. Therefore $G\cong{T}(1,0;h_1,h_2,h_3,h_4,h_5)$, where $h_i\ge 0$ for $1\le i\le 5$.
 
 \noindent\textbf{Subcase 1.3.} $s=1$ and $t=1$.
 
   In this case, since $F_{4}$, $F_{5}$ and $F_{6}$ are forbidden subgraphs, we obtain that there is a pendant path at the vertex of degree two in the three triangles of $T_1(1,1)$. Therefore $G\cong{T}(1,1;h_1,h_2,h_3,h_4,h_5)$, where $h_i\ge 0$ for $1\le i\le 5$.
   
\noindent\textbf{Conclusion 1.} When $G$ is a based on $T_1(s,t)$ graph, we obtain $G\cong{T}(s,t;h_1,,h_2,h_3,h_4,h_5)$ or $G\cong T_5$, where $t\ge 0 $, $s \ge 0 $ and $h_i\ge 0$ for $1\le i\le 5$. Given that ${T}(s,t;h_1,h_2,h_3,h_4,h_5)$ is a loose block graph, by the second conclusion of the Lemma \ref{th:ch-1}, we obtain
${T}(s,t;h_1,h_2,h_3,h_4,h_5)$ is  a connected  tricyclic graph with the  second  largest  distance eigenvalue less than \(-\frac{1}{2}\).

\noindent\textbf{Case 2.}  When $G$ is a based on $T_2(p,q)$ graph, we discuss the following three subcases.\\
\textbf{Subcase 2.1.} $p=0$ and $q=0$.\

  In this subcase, we obtain $G\cong T_3^k$ with $k\ge 0$, since $F_7$ is a forbidden subgarph.
  
  Next, we will prove graph $T_3^k$ with $k\ge 0$ is a tricyclic graph whose
  second largest distance eigenvalue less than \(-\frac{1}{2}\). By direct calculation, we have verified  that $-3$ is a distance eigenvalue of the graph $T_3^0$ with multiplicity 2, note that $T_3^0$ is a distance-preserving induced subgarph of $T_3^k$, according to the Lemma \ref{th:ch-2}, we obtain $-3$ is a distance eigenvalue of graph $T_3^k$.
   
   Assume that $G=T_3^k$ with $k\ge 1 $($n=|V(G)|=k+7$), and let $T_{3}^{k}$ be the tricyclic graph construted from the three triangles $\varDelta_{31}=v v_1v_2$, $\varDelta_{32}=v u_1u_2$  $\varDelta_{33}=vw_1w_2$, and attaching $k$ pendant edges at the vertex $v$ with maximum degree. we partition $V(G)=\{v\}\cup \{v_1,u_1,w_1\}\cup \{v_2,u_2,w_2\}\cup V_1$, where $V_1$ is the set of vertices with degree one. Based on this partition, we obtain
\[
D(G)+2I_n =
\begin{pmatrix}
	2 & J_{1\times 3} & J_{1\times 3}  & J_{1\times k} \\
	J_{3\times 1}  &  2J_{3\times 3}  & 2J_{3\times 3}-I_3 & 2J_{3\times k}\\
	J_{3\times 1} & 2J_{3\times 3}-I_3 & 2J_{3\times 3}  & 2J_{3\times k} \\
	J_{k\times 1} & 2J_{k\times 3} & 2J_{k\times 3}  & 2J_{k \times k} \\
\end{pmatrix}
\]

\[
D(G)+I_n =
\begin{pmatrix}
	1 & J_{1\times 3} & J_{1\times 3}  & J_{1\times k} \\
J_{3\times 1}  &  2J_{3\times 3}-I_3  & 2J_{3\times 3}-I_3 & 2J_{3\times k}\\
J_{3\times 1} & 2J_{3\times 3}-I_3 & 2J_{3\times 3}-I_3  & 2J_{3\times k} \\
J_{k\times 1} & 2J_{k\times 3} & 2J_{k\times 3}  & 2J_{k \times k}-I_k \\
\end{pmatrix}
\]

Note that the rank of $D(G)+2I_n$ is $8$, then $-2$ is a eigenvalue of $D(G)$ with multiplicity $k-1$. Analously, $-1$ is a eigenvalue of $D(G)$ with multiplicity $3$. Given that the above partition for $D(G)$ is equitable, by the Lemma \ref{th:ch-3}, then the eigenvalues of its quotient matrix $Q$ are also the eigenvalues of $D(G)$. $Q$ is shown below.
\[
Q =
\begin{pmatrix}
	0 & 3 & 3  & k \\
	1 &  4  & 5 & 2k\\
	1 & 5 & 4  & 2k \\
    1 & 6 & 6  & 2k-2 \\
\end{pmatrix}
\]
Define $f(\lambda)=det(\lambda I_4-Q)$, with direct calculation, we obtain
$$f(\lambda )=\lambda^4 -(2k + 6)\lambda^3 -(9k+31)\lambda^2 -(10k+36)\lambda - 3k - 12$$
 Note that $f(+\infty)>0$, $f(1)=-24k-84<0$, $f(-\frac{1}{2})=-\frac{15}{16}<0$, $f(-\frac{3}{4})=\frac{9k}{32}+\frac{105}{256}>0$, $f(-\frac{3}{2})=-\frac{3k}{2}-\frac{39}{16}<0$, $f(-4)=21k+276>0$. According to the above calculation, we obtain that the second largest root $ \lambda ^{\left( 2 \right)}$ of $f(\lambda)=0$ satisfies $-\frac{3}{4}<\lambda ^{\left( 2 \right)}<-\frac{1}{2}$. Based on the above argument, $ \lambda ^{\left( 2 \right)}$ is the second largest eigenvalues of $D(G)$. Therefore $\lambda _2(G)<-\frac{1}{2}$.

\noindent\textbf{Subcase 2.2.} $p\geq1$ or $q\geq1$.\

Given that $F_7$ and $F_{3}$ are forbidden subgarphs, then no graph which satisfies  \( \lambda_2(G) \) $<$ \(-\frac{1}{2}\) in this case.

\noindent\textbf{Conclusion 2.} When $G$ is a based on ${T}_2(p,q)$ graph, we obtain $G\cong T_3^k$ with $k\ge 0 $.

\noindent\textbf{Case 3.} When $G$ is a based on  $\text{T}_{4}^{t} $ graph, we discuss the following two subcases. Assume that $G=T_4^t$ with  $t\ge 0 $ ($n=|V(G)|=t+6$), let $T_{4}^{t}$ be the tricyclic graph obtained from the three triangles $\varDelta_{41}=v v_1v_2$, $\varDelta_{42}=v u_1u_2$, $\varDelta_{43}=vu_2u_3$, and attaching $t$ pendant edges at the vertex $v$ with maximum degree.\\
\textbf{Subcase 3.1.} $t=0$.\
 
 Because $F_9$ and $F_{10}$ are  forbidden subgarphs, we obtain $d(u_1)=d(u_3)=d(v_1)=d(v_2)=2$. Given that $F_8$ is a forbidden subgarph,
we obtain $d(u_2)\le 4$.

Assume that $d(u_2)=4$, we obtain $d(v)\le 6$, since $F_1$ is a forbidden subgarph. When $d(v)=6$, let $N_G\left( v \right)\setminus \{u_1,u_2,u_3,v_1,v_2\}=x_{1}$, We obtain $d(x_1)=1$, because $F_{10}$ is a forbidden subgarph. Therefore $G \cong T_{6}$. When $d(v)=5$, we obtain $G \cong T_{7}$.

Assume that $d(u_2)=3$, we consider the next subcase 3.2.

\noindent\textbf{Subcase 3.2.} $t\geq1$.

According to the Lemma \ref{th:ch-4}, given that $T_{4}^{t} $ is a distance-preserving  subgraph of $T_{4}^{t+1}$, we obtain $\lambda _2(T_4^t)\le \lambda _2(T_4^{t+1})$, which indicates that the sequence $\{\lambda _2(T_4^t): t=0,1,...\}$ is increasing or constant. Therefore it is sufficient to prove that \( \lambda_2(T_4^t) \) $<$ \(-\frac{1}{2}\) for large enough $t$.

Under this case, we partition $V(G)=\{v\}\cup \{v_1,v_2\}\cup \{u_1,u_3\}\cup \{u_2\}\cup V_1$, where $V_1$ is the set of vertices with degree one. Under this partition, we obtain

 \[
 D(G)+2I_n =
 \begin{pmatrix}
 	2 & J_{1\times 2} & J_{1\times 2} & 1 & J_{1\times t} \\
 	J_{2\times 1}  & J_{2\times 2}+I_2  & 2J_{2\times 2} & 2J_{2\times 1} & 2J_{2\times t}\\
 	J_{2\times 1} & 2J_{2\times 2} & 2J_{2\times 2}  & J_{2\times 1} & 2J_{2\times t} \\
 	1 & 2J_{1\times 2} & J_{1\times 2} & 2 & 2J_{1\times t} \\
 	J_{t\times 1} & 2J_{t\times 2} & 2J_{t\times 2} & 2J_{t\times1} & 2J_{t\times t}
 \end{pmatrix}
 \]
 
 \[
 D(G)+I_n =
 \begin{pmatrix}
 	1 & J_{1\times 2} & J_{1\times 2} & 1 & J_{1\times t} \\
 	J_{2\times 1}  & J_{2\times 2} & 2J_{2\times 2} & 2J_{2\times 1} & 2J_{2\times t}\\
 	J_{2\times 1} & 2J_{2\times 2} & 2J_{2\times 2}-I_2  & J_{2\times 1} & 2J_{2\times t} \\
 	1 & 2J_{1\times 2} & J_{1\times 2} & 1 & 2J_{1\times t} \\
 	J_{t\times 1} & 2J_{t\times 2} & 2J_{t\times 2} & 2J_{t\times1} & 2J_{t\times t}-I_t
 \end{pmatrix}
 \]
 Obviously, the rank of $D(G)+2I_n$ is $6$, then $-2$ is a eigenvalue of $D(G)$ with multiplicity $t$. Analously, $-1$ is a eigenvalue of $D(G)$ with multiplicity $1$. Note that the above partition for $D(G)$ is equitable, according to the Lemma \ref{th:ch-3}, then the eigenvalues of its quotient matrix $Q$ are also the eigenvalues of $D(G)$. $Q$ is shown below.
 \[
 Q =
 \begin{pmatrix}
 	0 & 2 & 2 & 1 & t \\
   1  & 1 & 4 & 2 & 2t\\
   1 & 4 & 2 & 1 & 2t \\
  1 & 4 & 2 & 0 & 2t \\
 1 & 4 & 4 & 2 & 2t-2  \\
 \end{pmatrix}
 \]
 Define $f(\lambda)=det(\lambda I_5-Q)$, with direct calculation, we obtain
 $$f(\lambda )=\lambda^5 - (2t+1)\lambda^4 - (15t+35)\lambda^3 - (35t+91)\lambda^2 - (26t+76)\lambda - 6t - 20$$
 Observe that $f(+\infty)>0$, $f(0)=-6t-20<0$, $f(-\frac{1}{2})=-\frac{15}{32}<0$, $f(-0.55)=\frac{2009t}{8000}-\frac{147871}{3200000}>0$ ($t\geq2$), $f(-1)=-2t-2<0$, $f(-2)=-6t<0$, $f(-3)=10>0$, $f(-4)=-14t-212<0$. According to the above calculation, we obtain that the second largest root $ \lambda ^{\left( 2 \right)}$ of $f(\lambda)=0 $ satisfies $-0.55<\lambda ^{\left( 2 \right)}<-\frac{1}{2}$. Based on the above argument, $\lambda ^{\left( 2 \right)}$ is the second largest eigenvalues of $D(G)$. Therefore $\lambda _2(G)<-\frac{1}{2}$.
 
 \noindent\textbf{Conclusion 3.} When $G$ is a based on  $T_{4}^t$ graph, we obtain $G \cong T_{6}$ or $G \cong T_{7}$ $(t=0)$, and $G \cong \text{T}_{4}^{t} $ ($t\geq1$).
 
According to the above arguments, we have the results as desired.
\qed\end{Tproof}

 \noindent\textbf{Declaration of competing interest}
 
 There is no competing interest.


\begin{thebibliography}{99}
\bibitem{RF1}A. Alhevaz, M. Baghipur, H.A. Ganie, On the second largest eigenvalue of the generalized distance matrix of graphs, Linear Algebra Appl. 603 (2020) 226–241.
\bibitem{RF2}M. Aouchiche, P. Hansen, Distance spectra of graphs: a survey, Linear Algebra Appl. 458 (2014) 301–386.
\bibitem{RF3}A.E. Brouwer, W.H. Haemers, Spectra of Graphs, Springer, New York, 2012. 
\bibitem{RF4}S. Fajtlowicz, Written on the wall: conjectures derived on the basis of the program Galatea Gabriella Graffiti, Technical report, University of Houston, 1998. 
\bibitem{RF5}R.L. Graham, H.O. Pollak, On the addressing problem for loop switching, Bell Syst. Tech. J. 50 (1971) 2495–2519. 
\bibitem{RF6}H.Y. Guo and B. Zhou, Graphs for which the second largest distance eigenvalue is less than {$-\frac{1}{2}$}, Discrete Mathematics 347(9)(2024) 114082,11 pp. 
\bibitem{RF7}W.H. Haemers, Interlacing eigenvalues and graphs, Linear Algebra Appl. 226(228)(1995) 593–616.
\bibitem{RF8}H.Q. Lin, Proof of a conjecture involving the second largest $D$-eigenvalue and the number of triangles, Linear Algebra Appl. 472 (2015) 48–53.
\bibitem{RF9}R.F. Liu, J. Xue, L.T. Guo, On the second largest distance eigenvalue of a graph, Linear Multilinear Algebra 65 (2017) 1011–1021.
\bibitem{RF10}H.Q. L, J.L. Shu, J. Xue, Y.K. Zhang, A survey on distance spectra of graphs, Adv. Math. (China) 50(1)(2021) 29--76.
\bibitem{RF11}R. Merris, The distance spectrum of a tree, J. Graph Theory 14 (1990) 365–369.
\bibitem{RF12}R.D. Xing, B. Zhou, On the second largest distance eigenvalue, Linear Multilinear Algebra 64 (2016) 1887–1898.
\bibitem{RF13}R.D. Xing, B. Zhou, On the two largest distance eigenvalues of graph powers, Inf. Process. Lett. 119 (2017) 39–43.
\bibitem{RF14}J. Xue, H.Q. Lin, J.L. Shu, On the second largest distance eigenvalue of a block graph, Linear Algebra Appl. 591 (2020) 284–298.
\bibitem{RF15} Z.X. Zhu, L.X. Tan and Z.~Y. Qiu, Tricyclic graph with maximal Estrada index, Discrete Appl. Math. 162 (2014) 364--372.

\end{thebibliography}
\end{document}